\documentclass[12pt,leqno]{amsart}
\usepackage{amsmath,amsthm, amscd, amssymb, amsfonts, mathrsfs,amsaddr}
\usepackage[all]{xy}
\usepackage{color}
\usepackage{enumitem}
\usepackage{mathdots}

\usepackage{hhline}
\usepackage[active]{srcltx}

\newcommand{\PSL}{\mathbf{PSL}}
\newcommand{\Gb}{\boldsymbol{G}}

\newcommand{\PSU}{\mathbf{PSU}}

\newcommand{\SU}{\mathbf{SU}}

\newcommand{\Jf}{{\tt J}}

\newcommand{\Irr}{\operatorname{Irr}}

\newcommand{\Tb}{\boldsymbol{T}}

\newcommand\w{\dot{w}}

\newcommand\toba{\mathfrak B }

\newcommand{\trid}{\triangleright}

\newcommand{\G}{{\mathbb G}}
\newcommand{\B}{{\mathbb{B}}}
\newcommand{\T}{{\mathbb{T}}}
\newcommand{\U}{\mathbb{U}}

\newcommand{\F}{{\mathbb F}}

\newcommand{\Gbtre}{{^2\!G_2(3)}}

\newcommand{\GL}{\mathbf{GL}}

\newcommand{\SL}{\mathbf{SL}}
\newcommand{\Sp}{\mathbf{Sp}}

\newcommand{\Fr}{\operatorname{Fr}}

\newcommand{\Oc}{{\mathcal O}}

\newcommand{\Aut}{\operatorname{Aut}}

\newcommand\Ad{\operatorname{Ad}}

\newcommand\Hom{\operatorname{Hom}}

\numberwithin{equation}{section}

\theoremstyle{plain}

\newtheorem{lema}{Lemma}[section]
\newtheorem{theorem}[lema]{Theorem}

\newtheorem{conjecture}[lema]{Conjecture}
\newtheorem{prop}[lema]{Proposition}

\newtheorem{question-app}{Question}

\theoremstyle{definition}
\newtheorem{definition}[lema]{Definition}

\theoremstyle{remark}
\newtheorem{obs}[lema]{Remark}

\newcommand\id{\operatorname{id}}

\newcommand\s{\mathbb S}

\def\pf{\begin{proof}}
\def\epf{\end{proof}}

\theoremstyle{remark}

\newcounter{tabla}\stepcounter{tabla}

\begin{document}

\renewcommand{\baselinestretch}{1.2}

\thispagestyle{empty}
\thispagestyle{empty}
\vspace*{-.2in}
\title[Nichols algebras over Suzuki and Ree groups]
{Finite-dimensional pointed Hopf algebras\newline over 
finite simple groups of Lie type VI. \newline 
Suzuki and Ree groups}

\author[G. Carnovale, M. Costantini]
{Giovanna Carnovale and Mauro Costantini}

\thanks{
	This work was partially supported by Progetto BIRD179758/17 of 
	the University of Padova.}

\address[A1]{
	Dipartimento di Matematica Tullio Levi-Civita,\\
	Universit\`a degli Studi di Padova,\\
	via Trieste 63, 35121 Padova, Italia\\
carnoval@math.unipd.it}
\address[A2]{
	Dipartimento di Matematica Tullio Levi-Civita,\\
	Universit\`a degli Studi di Padova,\\
	via Trieste 63, 35121 Padova, Italia\\
costantini@math.unipd.it}

\begin{abstract}
We analyse the rack structure of conjugacy classes in simple Suzuki and Ree groups and determine which classes are kthulhu. Combining these results with abelian rack techniques, we show that the only finite-dimensional complex pointed Hopf algebras over the simple Suzuki and Ree groups are their group algebras.
\end{abstract}
\maketitle
\setcounter{tocdepth}{1}


\section{Introduction}
This paper is part of an ongoing project with N. Andruskiewitsch and G. A. Garc\'ia, aimed at understanding finite-dimensional complex pointed Hopf algebras whose group of grouplikes is a finite simple group of Lie type \cite{ACG-I,ACG-II,ACG-III,ACG-IV, ACG-V}. We adopt notation and terminology from these papers, and for further details the reader is referred to them. We recall that a finite-dimensional pointed Hopf algebra  ${\mathcal H}$ has a natural filtration whose associated graded contains a graded associative algebra, the so-called Nichols algebra $\toba(V)$, whose structure depends on a representation $V$ of the finite group $\Gb$ of grouplike elements of ${\mathcal H}$ and a compatible $\Gb$-grading on $V$ (i.e., $V$ is a Yetter-Drinfeld module of $\Gb$). It is therefore crucial for our purposes to classify finite-dimensional Nichols algebras for Yetter-Drinfeld modules of $\Gb$. 

We recall the following folklore conjecture
\begin{conjecture}\label{conj:simple}Let $\Gb$ be a finite simple non-abelian group. Then, $\dim\toba(V)=\infty$ for every complex Yetter-Drinfeld module $V$ of $\Gb$. Thus, the only finite-dimensional complex pointed Hopf algebra whose group of grouplikes is $\Gb$ is the group algebra ${\mathbb C}\Gb$. 
\end{conjecture}

In fact Nichols algebras can be defined in a more general setup, whenever we have a vector space $V$ and an endomorphism $c\in {\rm GL}(V^{\otimes2})$ satisfying the braid equation $(c\otimes \id)(\id\otimes c)(c\otimes \id)=(\id\otimes c)(c\otimes \id)(\id\otimes c)$ as a  suitable quotient of the tensor algebra $T(V)$. If $V$ is a Yetter-Drinfeld module of $\Gb$, then $c$ is defined by
\begin{align}\label{eq:braiding}
c(v\otimes v')=g\cdot v'\otimes v,\quad v\in V_g,\,v'\in V,\quad V=\bigoplus_{g\in G}V_g.
\end{align}
The presentation of $\toba(V)$ depends only on the braiding $c$ and not on the Yetter-Drinfeld module structure itself. Indeed, such a braiding can arise from different actions of different groups on $V$ and $c$ can be described in terms of a combinatorial object called rack and a cocycle for it, \cite{AG}. Since important properties of racks are preserved by rack inclusions and projections, an approach to Nichols algebras done rack-by-rack might be more convenient than an approach group-by-group. In particular, the reduction to a (simple) rack  is relevant for the problem of classifying finite-dimensional pointed Hopf algebras whose group of grouplikes is finite but not simple. 

In our situation, a simple rack will always be a conjugacy class $\Oc$ in $\Gb$, with rack structure $g\trid h=ghg^{-1}$ for $g,h\in\Oc$. Conjugacy classes in different groups can be isomorphic as racks (e.g. unipotent conjugacy classes  arising from isogenous algebraic groups). An important goal is to classify finite-dimensional Nichols algebras for every conjugacy class in $\Gb$ and \emph{every cocycle}, and not only those coming from a Yetter-Drinfeld module of $\Gb$.

A series of conditions on racks (called type D, F and C) ensuring that the associated Nichols algebra is infinite-dimensional for any choice of a cocycle were given in \cite{AFGV-ampa,ACG-I,ACG-III}. In this case we say that the rack \emph{collapses}. In group theoretic terms these conditions are easy to state and well-behaved when passing to subgroups and quotients. The conjugacy classes that are not of type C, D, or F are called kthulhu: they are essentially those for which the possible non-empty intersections with a subgroup $H\leq\Gb$ are either a single conjugacy class in $H$ or consist of a set of commuting elements. For these classes we have no general strategy to deal with all cocycles. Yet, one can use the techniques developed in \cite{AFGV-ampa,AFGV-sporadic,AZ,FGV1} and the classification in \cite{H} to deal with the Nichols algebras associated with a kthulhu class and a cocycle coming from a (simple) Yetter-Drinfeld module of $\Gb$.
These techniques are often enough for dealing with Hopf algebras over $\Gb$ but might not propagate when passing to overgroups or groups projecting on $\Gb$. 

This paper deals with conjugacy classes in simple Suzuki groups $^2\!B_2(q)$, $q=2^{2h+1}$, $h\geq1$  and Ree groups $^2\!F_4(q)$, $q=2^{2h+1}$, $h\geq1$ and $^2G_2(q)$, $q=3^{2h+1}$, $h\geq1$. They were firstly presented and studied in \cite{Su,ree1,ree2} but we will use the uniform description in terms of fixed points of a Steinberg endomorphism as in \cite{Carter1,MT}. Working group-by-group we prove the following result:
\begin{theorem}Let $\Gb$ be a simple Ree or Suzuki group and let $\Oc$ be the class of  $x\in\Gb$. Then $\Oc$ is kthulhu if and only if
\begin{align*}
&\Gb=\,^2\!B_2(q)\mbox{ and }|x|=2  \mbox{ or }|x|\mbox{ divides }q^2+1;\\
&\Gb=\,^2\!F_4(q)\mbox{ and }|x|\neq13\mbox{ and  divides }q^4-q^2+1;\\
&\Gb=\,^2\!G_2(q)\mbox{ and }|x|=3\mbox{ or }|x|\neq7\mbox{ and divides }q^2-q+1.\\
\end{align*}
\end{theorem}
In other words, kthulhu classes are either classes of unipotent elements of prime order, or they are represented by elements in maximal tori whose normaliser is the only maximal subgroup containing them.  Up to some exception of small order, every non-trivial element in such a torus lies in a kthulhu class. 
%
Then we focus on kthulhu classes in each group and finally we prove: 
\begin{theorem}\label{thm:main}Conjecture  \ref{conj:simple} holds if $\Gb$ is a  simple Suzuki or Ree group.
\end{theorem}
%
%

\section{Notation and background}
For any automorphism $\sigma$ of an algebraic structure $X$, we shall denote by $X^\sigma$ the set of elements fixed by $\sigma$. 
For $G$ a group, the orbit  of an element $g\in G$ under the conjugation action of a subgroup $H\leq G$ will be denoted by $\Oc_g^H$. The superscript will be omitted if the ambient group is clear from the context. To keep uniformity with the previous papers in the series, we will denote the conjugation action by: $g\trid h:=ghg^{-1}$. The centraliser of an element $x\in G$ will be denoted by $C_G(x)$, and the set of isomorphism classes of irreducible representations of a group $H$ will be denoted by $\Irr(H)$.

\subsection{Preliminaries on racks}\label{sec:racks}
In this section we introduce some preliminary notions on the rack structure specialised to the case of a conjugacy class.

\begin{definition}(\cite[Definition 2.3]{ACG-III}, \cite[Definition 3.5]{AFGV-sporadic}, \cite[Definition 2.4]{ACG-I}). A conjugacy class $\Oc$ in a finite group $M$ is said to be of type
\begin{enumerate}[leftmargin=*]
\item[\text{C}]  if there are $H\leq M$ and $r,\,s\in \Oc\cap H$ 
such that
\begin{enumerate}
\item $\Oc_r^H\neq \Oc_s^H$,
\item $rs\neq sr$,
\item $H=\langle \Oc_r^H,\, \Oc_s^H\rangle$,
\item either $\min(| \Oc_r^H|,\, |\Oc_s^H|)>2$ or $\max(|\Oc_r^H|,\, |\Oc_s^H|)>4$;
\end{enumerate}
\item[\text{D}] if there are $r,\,s\in\Oc$ such that
\begin{enumerate}
\item $\Oc_r^{\langle r,\,s\rangle}\neq \Oc_s^{\langle r,\,s\rangle}$,
\item $(rs)^2\neq(sr)^2$;
\end{enumerate}
\item[\text{F}]
 if there are $r_i\in\Oc$, for $1\leq i\leq4$ such that
\begin{enumerate}
\item $\Oc_{r_i}^{\langle r_i,\,1\leq i\leq 4\rangle}\neq \Oc_{r_j}^{\langle r_i,\,1\leq i\leq 4\rangle }$, for $i\neq j$,
\item $r_ir_j\neq r_jr_i$, for $i\neq j$.
\end{enumerate}
\end{enumerate}
A conjugacy class is called \textit{kthulhu} if it is of none of these types.
\end{definition}
The relevance of the above conditions relies on the following results, that we apply to the special case of conjugacy classes.

\begin{prop}\label{prop:type_Nichols}(\cite[\S 2.2]{ACG-III}, \cite[\S 3.2]{AFGV-sporadic},\cite[\S 2.2]{ACG-I}).
\begin{enumerate}
\item If a rack $X$ is of type C, D, or F, then $\dim\toba(X,\mathbf{q})=\infty$ for every cocycle $\mathbf{q}$ for $X$, i.e., it collapses. 
\item\label{item:inj-proj} If a rack contains or projects onto a rack of  type C, D, or F, then it is of the same type.
\end{enumerate}
\end{prop}

\begin{obs}\label{obs:conditions}
\begin{enumerate}
\item\label{item:andr} Assume that $M$ is a finite group with $M/Z(M)$ simple non-abelian. If $m\in M\setminus Z(M)$, then there exists $g\in M$ such that $[g\trid m,m]\neq 1$. Otherwise, $N:=\langle \Oc_{m}^{M}\rangle$ would be an abelian normal subgroup of $M$ and $N/(Z(M)\cap N)$ would be an abelian normal subgroup of $M/Z(M)$. Therefore $N$ would be central, while $m\not\in Z(M)$, a contradiction. 
\item\label{item:bigger} If $r\in M$ with $|r|$ odd, and $r,s\in \Oc_r^M$ satisfy $rs\neq sr$, then for any $H\leq M$ such that $\langle r,\,s\rangle \leq H$ we have $\min(| \Oc_r^H|,\, |\Oc_s^H|)>2$. Indeed, 
$3\leq |\Oc_s^{\langle r\rangle}|\leq | \Oc_s^H|$ and  $3\leq |\Oc_r^{\langle s\rangle}|\leq | \Oc_r^H|$. 
\end{enumerate}
\end{obs}

\begin{lema}\label{lem:product}Assume $M=M_1\times M_2$ is a finite group such that $M_1/Z(M_1)$ is simple non-abelian and let $m_i\in M_i\setminus Z(M_i)$ for $i=1,2$. If $|m_1|$ is odd, $|m_2|\neq 2$ and $m_2$ is real in $M_2$, then $\Oc_{(m_1,m_2)}^M$ is of type C.
\end{lema}
\pf By Remark \ref{obs:conditions} \eqref{item:andr} there is $g_1\in M_1$ such that $[g_1\trid m_1,m_1]\neq1$. Let $g_2\in M_2$ be such that $g_2\trid m_2=m_2^{-1}$. We set $r=(m_1,m_2)$, $s:=(g_1,g_2)\trid(m_1,m_2)$ and $H=\langle r,\,s\rangle\leq \langle m_1, g_1\trid m_1\rangle\times\langle m_2\rangle$. By construction $rs\neq sr$ and $H=\langle \Oc_r^H,\Oc_s^H\rangle$. The inequality $m_2^2\neq 1$ implies $\Oc_r^H\cap\Oc_s^H=\emptyset$. In addition $|\Oc_r^H|=|\Oc_{m_1}^{\langle m_1, g_1\trid m_1\rangle}|\geq|\Oc_{m_1}^{\langle g_1\trid m_1\rangle}|\geq 3$ because $|m_1|=|g_1\trid m_1|\geq 3$, and similarly for $|\Oc_s^H|$.
\epf

In the following Remark we recall  an argument used in \cite[Proposition 5.5]{ACG-III} in order to prove that certain classes are kthulhu. 
\begin{obs}\label{obs:inductive_kthulhu}Let $\Oc$ be a conjugacy class in a finite group $G$. Assume that for any $H\leq G$ the intersection $\Oc\cap H$ is either empty, a unique conjugacy class in $H$, or consists of mutually commuting elements. Then $\Oc$ is kthulhu. We usually deal with intersections with subgroups using the list of maximal subgroups as follows. For any maximal  $M< G$  we analyse $\Oc\cap M$. If $\Oc\cap M=\emptyset$ or consists of commuting elements, it will be again so for any $H\leq M$. Then we show that for the remaining subgroups $\Oc\cap M$ is a single conjugacy class in $M$ and observe that in this case the structure of $M$ and of its (maximal) subgroups are well-understood. In most cases $M$ will be a finite simple group of Lie type of the same sort as $G$ but over a smaller field, or $\PSL_2(q)$. This way we reduce from the pair $(\Oc, G)$ to the pair $(\Oc\cap M, M)$. The class $\Oc\cap M$ will usually have the same features as $\Oc$ had in $G$ and we proceed inductively.
\end{obs}

In order to implement the above mentioned analysis, we will make use of the following standard observation.

\begin{obs}\label{obs:split} Let $M=N\rtimes\langle a\rangle$ be a  finite group with $C_N(a)=\{1\}$.
\begin{enumerate}
\item\label{item:mauro1}We have the equality $\Oc_{a}^M=Na$. Indeed, the inclusion $\subset$  follows from normality of $N$, and $C_N(a)=\{1\}$ implies that the two sets have the same cardinality. 
\item\label{item:mauro2} For any $H\leq M$ such that $a\in H$, then $\Oc_a^H=(H\cap N)a$. Indeed, if $a\in H$, then $H=(N\cap H)\rtimes\langle a\rangle$, so \eqref{item:mauro1} applies.
\item\label{item:mauro3} If $G$ is a finite group containing $M$ and such that $\Oc_a^G\cap M\subset Na$, then for any $H\leq M$ with $a\in H$ we have $\Oc_a^G\cap H=\Oc_a^H$.
Indeed,  by \eqref{item:mauro2} there holds $\Oc_a^G\cap H\subset (N\cap H)a=\Oc_a^H\subset \Oc_a^G\cap H$, whence the equality.
\end{enumerate}
\end{obs}

\subsection{Nichols algebras of Yetter-Drinfeld modules and abelian subracks}\label{sec:abelian}

In this section we provide some necessary ingredients for dealing with Conjecture \ref{conj:simple}.  The first key observation is the following. 
\begin{obs}\label{obs:reduction} (\cite[\S 1.2]{AFGV-sporadic}) Let $H$ be a finite group. If  $\dim\toba(V)=\infty$ for every simple  Yetter-Drinfeld module $V$ of $H$, then $\dim\toba(V')=\infty$ for every  Yetter-Drinfeld module $V'$ of $H$.
\end{obs}

Simple Yetter-Drinfeld modules of $H$ are parametrized by pairs $(\Oc,\rho)$ where the $\Oc=\Oc_g^H$ is a conjugacy class in $H$ and is the support of the grading, and $\rho\in{\rm Irr}(C_H(g))$. If $\rho\colon C_H(g)\to \GL(W)$, 
the corresponding simple Yetter-Drinfeld module $M(\Oc,\rho)$ has $H$-module structure and grading defined by:
\begin{align}\label{eq:YD}
 M(\Oc,\rho)={\rm Ind}_{C_H(g)}^H\rho={\mathbb C}H\otimes_{{\mathbb C}C_H(g)} W,&&m\otimes W\subset M(\Oc,\rho)_{m\trid g},\quad m\in H.
\end{align}
If $\Oc$ is of type C, D, or F, then Proposition \ref{prop:type_Nichols} (1) ensures that $\dim\toba(M(\Oc,\rho))=\infty$ for any choice of $\rho\in {\rm Irr}(C_H(g))$.  For a kthulhu conjugacy class $\Oc=\Oc_h^H$ the conclusion of Proposition \ref{prop:type_Nichols} (1) cannot be inferred. We recall a strategy developed in \cite{AFGV-ampa,AFGV-sporadic,AZ,FGV1} and references therein to estimate the dimension of $\toba(M(\Oc,\rho))$. 

Assume $A\leq C_H(g)$ is an abelian subgroup containing $g$. 
Then, $\Oc\cap A$ is an abelian subrack of $\Oc$ and 
$\rho(A)$ stabilises a line ${\mathbb C}v$ in $W$: we call $\chi$ its character. We set
\begin{align}\label{eq:braiding-abelian}
 \Oc\cap A=\{x_0:=g_0\trid g,\,\ldots, x_r:=g_r\trid g\},&&q_{ij}:=\chi(g_j^{-1}g_i\trid g).
\end{align}
Then $M_A={\rm span}_{\mathbb C}(g_i\otimes v, i=0,\ldots,r)$ is a braided subspace of $M(\Oc,\rho)$, i.e., $c(M_A\otimes M_A)=M_A\otimes M_A$, where $c$ is as in \eqref{eq:braiding}. The restriction of $c$ to $M_A\otimes M_A$ is given by $c((g_i\otimes v)\otimes (g_j\otimes v))=q_{ij}(g_j\otimes v)\otimes (g_i\otimes v)$, i.e., it is of diagonal type with $q_{ii}=\chi(g)$ for every $i$. Now $\toba(M_A)$ is a subalgebra of $\toba(M(\Oc,\rho))$. We can invoke the classification results for finite-dimensional Nichols algebras for braided spaces of diagonal type in \cite{H}, and if $\dim\toba(M_A)=\infty$ we can conclude that $\dim\toba(M(\Oc,\rho))=\infty$.

\medskip

Proposition \ref{prop:type_Nichols} implies that for an extension $G$ of a group $H$, a conjugacy class $\Oc$ in $G$ and a cocycle  $\mathbf{q}$ for $\Oc$, the dimension of $\toba(M(\Oc,\mathbf{q}))$ can be finite  only if the projection of $X$ in $H$ is kthulhu, giving indications for the quest of finite-dimensional Nichols algebras over $G$. However, information concerning $H$ needs to be integrated with other tools to retrieve complete information on $G$: on the one hand the lift of a kthulhu conjugacy class is not necessarily kthulhu, so the list of classes in $G$ potentially yielding finite-dimensional Nicholas algebras might be reduced; on the other hand, the techniques described in this section do not immediately lift from $H$ to $G$: for instance, the lift of $A$ may fail to be abelian.

\subsection{Construction of the groups}\label{subsec:construction}

Let $p$ be a prime, $h\geq0$, $q=p^{2h+1}$ and $\G$ a simply-connected simple algebraic group over $\overline{\F_p}$. We recall the construction of the groups $^{2}\!B_2(q)$, $^{2}\!F_4(q)$
and $^{2}\!G_2(q)$ from \cite[\S 13]{Carter1}, 
as fixed point sets of certain Steinberg endomorphism $F$ in $\G$. Let $\T$ be a fixed maximal torus in $\G$, with corresponding root system $\Phi$, root subgroups $\U_\alpha$ for $\alpha\in \Phi$, and Weyl group $W=N_{\G}(\T)/\T$. We fix an isomorphism $x_\alpha\colon \overline{\F_p}\to \U_\alpha$ for each $\alpha\in \Phi$ and a set of simple roots $\Delta$, with corresponding positive roots $\Phi^+$. The group $W$ acts by isometries on $E={\mathbb R}\otimes_{\mathbb Z}{\mathbb Z}\Phi$. 

We will focus on the cases in which the pair $(\Phi,p)$ is either $(B_2,2)$, $(F_4,2)$ or $(G_2,3)$. In the latter case we assume $x_\alpha$ is as in \cite[\S 12.4]{Carter1}.

The non-trivial symmetry of the Coxeter graph of $\G$ induces a permutation $\theta\colon \Phi\to \Phi$, \cite[\S 12.3, 12.4]{Carter1}. 
We denote by $\tau$ the unique involutory isometry of $E$ such that $\tau(\alpha)\in{\mathbb R}_{>0}\theta{\alpha}$ for all  $\alpha\in \Phi$:
\begin{align*}
\tau(\alpha)=\begin{cases}
\frac1{\sqrt2}\ \theta \alpha &\textrm{if $\alpha$ is short,}\\
\sqrt 2\ \theta \alpha&\textrm{if $\alpha$ is long,} 
\end{cases}\qquad
\textrm{for $\Phi$ of type $B_2$ or $F_4$},
\end{align*}
\begin{align*}
\tau(\alpha)=\begin{cases}
\frac1{\sqrt3}\ \theta \alpha &\textrm{if $\alpha$ is short,}\\
\sqrt 3\ \theta \alpha&\textrm{if $\alpha$ is long,} 
\end{cases}\qquad
\textrm{for $\Phi$ of type $G_2$}. 
\end{align*}

There is a graph automorphism $\vartheta$ of $\G$ preserving $\T$ and such that 
$\vartheta(\U_\alpha)=\U_{\theta\alpha}$ for all $\alpha\in \Phi$, \cite[\S 12.3, 12.4]{Carter1}. It is defined as follows on root subgroups:
\begin{align*}
\vartheta(x_\alpha(\xi))&=\begin{cases}
x_{\theta \alpha}(\xi^p) &\qquad \textrm{if $\alpha$ is short,}\\
x_{\theta \alpha}(\xi) &\qquad\textrm{if $\alpha$ is long.}
\end{cases}&\xi\in \overline{\F_p}
\end{align*}

Let $\Fr_{p^h}$ be the field automorphism of $\G$ induced by the automorphism
$\lambda\mapsto\lambda^{p^h}$ of $\overline{\F_p}$ and let $F:\G\to\G$ be the Steinberg endomorphism $F=\vartheta\circ \Fr_{p^h}=\Fr_{p^h}\circ \vartheta$. Then $\Gb:=\G^{F}=\{x\in \G\mid F(x)=x\}$ are the Suzuki groups for $\G$ of type $B_2$ and the Ree groups for $\G$ of type $F_4$ or $G_2$.

Note that
\begin{align*}
F^2:x_\alpha(\xi)\mapsto x_{\alpha}(\xi^{q})
\end{align*}
for every $\alpha\in \Phi$, so that $\Gb$ is contained in 
$B_2({2^{2h+1}}),\ F_4({2^{2h+1}}),\ G_2({3^{2h+1}})$
respectively. For convenience, we denote $\Gb$ by
${}^2\!B_2(q),\ {}^2\!F_4(q),\ {}^2\!G_2(q)$
respectively. Notice that this is not the only notation for the Suzuki or Ree groups, often the convention $q^2=p^{2h+1}$ is used. We have
\begin{align*}
\left\vert {}^2\!B_2(q)\right\vert&=q^2(q-1)(q^2+1),&\left\vert {}^2\!G_2(q)\right\vert&=q^3(q-1)(q^3+1),\\
\left\vert {}^2\!F_4(q)\right\vert&=q^{12}(q-1)(q^3+1)(q^4-1)(q^6+1).
\end{align*}
We recall that these groups are simple for $h\geq1$.

Let $\B\leq\G$ be the Borel subgroup generated by $\T$ and the $\U_\alpha$, $\alpha\in\Phi^+$, let $\B^-\leq\G$ be the opposite Borel subgroup and let $\U$ and $\U^-$ be their unipotent radicals. Every unipotent conjugacy class in $\G$ is represented by an element in $\U$ and, for any fixed ordering of $\Phi^+$, every element in $\U$ can be uniquely written as a product $\prod_{\gamma\in\Phi^+}x_\gamma(c_\gamma)$ for $c_\gamma\in\overline{\F_q}$. Also $U:=\U^F$ and $U^-:=(\U^-)^F$ are Sylow $p$-subgroups of $\Gb$, \cite[Corollary 24.11]{MT} so all unipotent conjugacy classes (i.e. consisting of elements of order a power of $p$) in $\Gb$ intersect $U$ and $U^-$. We recall that every $F$-stable maximal torus in $\G$ is of the form $g\T g^{-1}$ for some $g\in \G$ such that $\w:=g^{-1}F(g)\in N(\T)$. Two such tori are $\Gb$-conjugate if and only if the corresponding Weyl group elements are $F$-conjugate, \cite[Proposition 25.1]{MT}. 
We denote by $\T_w$ a maximal torus whose associated Weyl group element is $w=\w\T\in W$ and we set $\Tb_w=\T_w^F$. Every semisimple element   (i.e. of order coprime with $p$) in $\Gb$ is contained in some $\Tb_w$, for some $w\in W$, \cite[Proposition 26.6]{MT}.
There is a formula for the order of $\Tb_w$. Let $Y=\Hom(\overline {\F_p}^\times,\T)$ be  the cocharacter group of $\T$. Then $W$ and $F$ act on $Y$, hence on $Y\otimes\mathbb R$ and 
$\left\vert \Tb_w\right\vert =
\vert\textstyle\det_{_{Y\otimes\mathbb R}}(w^{-1}\circ F-1)\vert$, see \cite[Proposition 25.3 (c)]{MT}.
An element $\sigma\in W$ has a representative in $N_{\Gb}(\T_w)$ if and only if $\sigma\in C_W(\tau w)$, and $|N_{\Gb}(\T_w)/\Tb_w|=|C_W(\tau w)|$,  \cite[Proposition 25.3 (a)]{MT}. 

\begin{obs}\label{obs:mixed}When dealing with mixed classes, i.e., classes  of elements $x\in\Gb$ that are neither semisimple nor unipotent we adopt the strategy developed in \cite[\S 3]{ACG-V}. Let $x=x_sx_u$ be the Jordan decomposition of $x$. We recall that $[C_{\G}(x_s),C_{\G}(x_s)]$ is a semisimple group whose root system has a base that can be indexed by a set of nodes $\Sigma$ in the extended Dynkin diagram of $\G$, \cite[Remark 14.5]{MT}. In addition, $\Sigma$ must be stable by ${\Ad}(\w)\circ F$ for some $w\in W$. Since $W$ preserves the root lengths and $\vartheta$  does not, if $\Sigma$ is non-empty, identified with the corresponding subset of $\Phi$, can only have the same amount of short and long roots, providing a strong restriction on the possibilities for $\Sigma$. Also, $[C_{\G}(x_s),C_{\G}(x_s)]^F\simeq \langle \T,\U_{\pm\alpha}~|~\alpha\in \Sigma\rangle^{\Ad(\w)F}$ and the following natural rack inclusion 
\begin{align}\label{eq:inclusion}
\Oc_{x_u}^{[C_{\G}(x_s),C_{\G}(x_s)]^F}\simeq x_s\Oc_{x_u}^{[C_{\G}(x_s),C_{\G}(x_s)]^F}=\Oc_{x_sx_u}^{[C_{\G}(x_s),C_{\G}(x_s)]^F}
\subset \Oc_{x_sx_u}^{\G^F} \end{align} shows that if $\Oc_{x_u}^{[C_{\G}(x_s),C_{\G}(x_s)]^F}$ is not kthulhu, then $\Oc_x^{\Gb}$ is again so.
\end{obs}

\begin{obs}\label{obs:ss-real}
Since in $B_2,\,F_4$ and $G_2$ the longest element $w_0$ in $W$ is $-\id$, for any $w\in W$ there  is always a representative of $w_0$ in $N_{\Gb}(\T_w)$ and therefore all semisimple classes in $\Gb$ are real. 
\end{obs}


\begin{obs}\label{obs:order-tori}
Since $X\pm1$ divides $X^{2m+1}
\pm 1$,
if  $q=q_0^{2m+1}$ and $(d,q^{k}\pm1)=1$, then $(d, q_0^k\pm1)=1$. So, if an  element of $^{2}\!B_2(q)$, $^{2}\!F_4(q)$ or $^{2}\!G_2(q)$ lies in a torus whose order is coprime to $q^k\pm1$, and it also lies in a subgroup isomorphic to $^{2}\!B_2(q_0)$, $^{2}\!F_4(q_0)$ or $^{2}\!G_2(q_0)$, respectively, then it will lie in a torus therein whose order is coprime to $q_0^k\pm1$.
\end{obs}

\begin{obs}\label{obs:ss-one-class}
If $F'$ is a Steinberg automorphism of $\G$ such that $\G^{F'}\leq \Gb$ and $g\in \Gb$ is semisimple, then 
 $\Oc_g^{\Gb}\cap \G^{F'}\subset \Oc_g^{\G}\cap\G^{F'}$. Since $\Oc_g^{\G}$ is semisimple, the right hand side is either empty or a unique semisimple conjugacy class in $\G^{F'}$ by \cite[\S 3.4 (c), p. 177]{sp-st}, so the same holds for the left hand side.
\end{obs}

\section{The Suzuki groups $^{2}\!B_2(2^{2h+1})$}

In this section $p=2$, $q=2^{2h+1}$, $h\geq0$, $\Gb={}^2\!B_2(q)$ a Suzuki group. We recall some basic facts  from \cite{Su}.

We will need the automorphism $\delta=\Fr_{2^{h+1}}$ of $\F_{q}$, so that $\delta^2=\Fr_2$ and $\F_{q}^\delta=\F_2$.  
\begin{obs}\label{obs:kappa}
\begin{enumerate}
\item 
Let $k\in\F_q$ be such that $k\delta(k)=1$. Then, $1=\delta(k\delta(k))=\delta(k)k^2$, so $k=k^2\in\F_2$.
\item The group morphism $\varphi\colon \F_q^\times \to \F_q^\times$ given by $\varphi(k)=k\delta(k)$ is injective, therefore it is an isomorphism. 
\end{enumerate}
\end{obs}

We realize $\G=\Sp_{4}(\overline{\F_2})$ as the group of matrices in $\GL_4(\overline{\F_2})$ preserving the bilinear form associated with the matrix
$\Jf=\left(\begin{smallmatrix}
0&0&0&1\\
0&0&1&0\\
0&1&0&0\\
1&0&0&0
\end{smallmatrix}\right)$. Then $\T$ can be chosen to be the subgroup of diagonal matrices and $\B$ the subgroup of upper-triangular matrices. The Sylow $2$-subgroup $U^-$ of $\Gb$ is given by the matrices of the form
\begin{align*}
U(a,\,b):=\left(\begin{matrix}1\\
a&1\\
a\delta(a)+b&\delta(a)&1\\
a^2\delta(a)+ab+\delta(b)&b&a&1\end{matrix}\right),&&a,b\in\F_q.
\end{align*}
with multiplication rule
\begin{align}\label{eq:productU}
U(a,b)U(c,d)=U(a+c,a\delta(c)+b+d),&&a,b,c,d\in\F_q,
\end{align}

For any $k\in\F_q^\times$ we have $t_k:={\rm diag}(\xi_1,\xi_2,\xi_2^{-1},\xi_1^{-1})\in\Tb$ where $\delta(\xi_1)=k\delta(k)$ and $\delta(\xi_2)=k$.
There holds:
\begin{align}\label{eq:comm_tU}
t_k^{-1}U(a,b)t_k=U(ak,bk\delta(k)),&&a,b\in\F_q,\,k\in\F_q^\times.
\end{align}

It follows from  \cite[Proposition 1, 2, 3, 7]{Su}  that
if $x$ is a non-trivial element of a Sylow $2$-subgroup $Q$ of $\Gb$, then $C_{\Gb}(x)\leq Q$. In particular, the order of the elements in $\Gb$ is either a power of $2$ or odd.

It follows from \eqref{eq:productU} that all non-trivial involutions are conjugate to an element of the form $U(0,b)$, so by Remark \ref{obs:kappa} and \eqref{eq:comm_tU} all non-trivial involutions are conjugate. 

The elements of odd order are semisimple, and therefore their conjugacy classes are represented by an element in a maximal torus $\Tb_w$, where $w$ runs through a set of representatives of $F$-conjugacy classes in $W$. Up to conjugacy they are: 
$\Tb$, of order $q-1$  and the two cyclic subgroups $\Tb_{s_1}$ and $\Tb_{s_1s_2s_1}$ whose orders are  $2^{2h+1}\pm2^{h+1}+1=q\pm\sqrt{2q}+1$, so $|\Tb|,  |\Tb_{s_1}|$ and $|\Tb_{s_1s_2s_1}|$ are mutually coprime. 

The maximal subgroups of $\Gb$ are the conjugates of the following subgroups, \cite[Theorems 9, 10]{Su}: 
\begin{enumerate}
 \item $B^-=\Tb\ltimes U^-$ of order $q^2(q-1)$;
 \item $N_{\Gb}(\T)$ of order $2(q-1)$;
 \item $N_{\Gb}(\T_{s_1})$ and $N_{\Gb}(\T_{s_1s_2s_1})$ of order 
$4(q\pm \sqrt{2q}+1)$;
\item $^2\!B_2(2^{2m+1})$ for $(2h+1)/(2m+1)$ a prime number. 
\end{enumerate}

%
%

\subsection{Collapsing racks}


\begin{lema}\label{lem:invo_suzuki}If $\Oc$ consists of elements of order $2$, then $\Oc$ is kthulhu.
\end{lema}
\pf The class $\Oc$ is contained in a unipotent class in  $\Sp_{4}(\F_q)$  and all classes of non-trivial involutions therein are kthulhu \cite[Lemma 4.22(2), Lemma 4.26]{ACG-II}, \cite[Lemma 2.14]{ACG-III}. We conclude by \cite[Lemma 2.5]{ACG-IV}.
\epf

\begin{lema}\label{lem:suz_order4}Assume $h>0$. If $\Oc$ consists elements of order $4$, then it is of type F.
\end{lema}
\pf By \eqref{eq:productU} any element of order $4$ has a representative $r=U(a,b)$ for some $a\neq0$, $a,b\in\F_q$. 
Since $h>0$, there are  distinct $k_j\in\F_q^\times$ for $j=0,1,2,3$ and we set $r_j:=t_{k_j}\trid r=U(ak_j,bk_j\delta(k_j))$. 
For any $c,d\in\F_q$ we have $U(c,d)^{-1}=U(c,d')$ for some $d'\in\F_q$ and
$U(c,d)\trid U(a,b)=U(a,b')$ for some $b'\in\F_q$. As $\langle r_i,\,i=0,1,2,3\rangle\leq U^-$, we deduce that
$\Oc_{r_i}^{\langle r_i,\,i=0,1,2,3\rangle}\neq \Oc_{r_j}^{\langle r_i,\,i=0,1,2,3\rangle}$ for $i\neq j$. 
In addition,
\begin{align*}
r_ir_j=U(a(k_i+k_j),a\delta(a)k_i\delta(k_j)+bk_i\delta(k_i)+bk_j\delta(k_j))
\end{align*}
so $r_ir_j=r_jr_i$ if and only if $k_ik_j^{-1}=\delta(k_ik_j^{-1})$ if and only if $k_i=k_j$ if and only if $i=j$.
Whence, $\Oc$ is of type F.
\epf

\begin{lema}\label{lem:toro-split} Let $1\neq t\in\Tb$.
Then $\Oc_{t}^{\Gb}$ is of type C. 
\end{lema}
\pf Recall that $t\neq t^{-1}\in \Oc_t^{\Gb}$ by Remark \ref{obs:ss-real} and that $C_{\Gb}(U(1,0))\leq U^-$,  \cite[Proposition 1, 2, 3, 7]{Su}. 
It follows from Remark \ref{obs:split} that $U^-\trid t=tU^-\subset \Oc_{t}^{\Gb}$ and similarly, $U^-\trid t^{-1}=t^{-1}U^-\subset \Oc_{t}^{\Gb}$. Let $H:=\langle t,U^-\rangle$, $r=t$, $s=t^{-1}U(1,0)\in\Oc_{t}^{\Gb}\cap H$. Then, $rs\neq sr$ by \eqref{eq:comm_tU}. Also, 
$\Oc_{r}^H=\Oc_{t}^{U^-}=tU^-$ and $\Oc_{s}^H=\Oc_{t^{-1}}^H=\Oc_{t^{-1}}^{U^-}=t^{-1}U^-$, the two sets are clearly disjoint and
\begin{align*}\langle \Oc_{r}^H,\Oc_{s}^H\rangle=\langle tU^-,t^{-1}U^-\rangle=\langle t, U^-\rangle=H.\end{align*}
In addition $| \Oc_{r}^H|=|\Oc_{s}^H|=|U^-|=q^2>2$, so $\Oc_{t}^{\Gb}$ is f type C.
\epf

\begin{obs}\label{obs:suz2}All classes in $^2\!B_2(2)$ are khtulhu. Indeed, $|^2\!B_2(2)|=20$ and its elements  have either order $2,4$ or $5$. We realize it  as the subgroup of matrices
$m(a,x)=\left(\begin{smallmatrix}
1&0\\
x&a\end{smallmatrix}\right)$ where $x\in\F_5$ and $a\in\F_5^\times$. 
Let $g\in \,^2\!B_2(2)$ with $|g|=5$. Then, $g=m(0,x)$ for some $x\in\F_5$ and all such elements form an abelian normal subgroup of $^2\!B_2(2)$. Hence, all elements in the class of $g$ commute, so classes of elements of order $5$ cannot be of type C, D nor F.  Let now $|g|=4$. Then $g=m(a,x)$, for $a=2$ or $3$ in $\F_5$ and some $x\in\F_5$ and $h\in\Oc_g$ if and only if $h=m(a,y)$ for some $y\in\F_5$. 
But then $\langle g,h\rangle$ contains elements of order $4$ and it is either cyclic of order $4$ or the whole group. Hence, $\Oc_g$ is kthulhu.
For non-trivial involutions we invoke Lemma \ref{lem:invo_suzuki}. 
\end{obs}

\begin{lema}\label{lem:suz-semis-kthulhu}The non-trivial classes represented by elements in the subgroups $\Tb_{s_1}$ and  $\Tb_{s_1s_2s_1}$ are kthulhu. \end{lema}
\pf We use the inductive argument from Remark \ref{obs:inductive_kthulhu}. When $h=0$, all classes are kthulhu by Remark \ref{obs:suz2}. So assume $h > 0$.
Let $\Tb'$ be one of these tori, let $g\in\Tb'\setminus\{1\}$ and let $\Oc=\Oc_g^{\Gb}$. Since $|g|$ divides $(q+ \sqrt{2q}+1)(q-\sqrt{2q}+1)=q^2+1$, a maximal subgroup $M$ of $\Gb$ meeting $\Oc$ cannot be conjugate to $B^-$ or $N_{\Gb}(\T)$. Also by order reasons, if it meets $N_{\Gb}(\T_{w})$ for $w=s_1$ or $s_1s_2s_1$, then $\Tb'=\Tb_w$ and since  $N_{\Gb}(\T_{w})\simeq \Tb_w\rtimes C_4$, all its elements of odd order are contained in $\Tb_w$. Thus, $\Oc\cap N_{\Gb}(\T_w)$ consists of commuting elements. We finally consider $M={}^2\!B_2(2^{2f+1})$ with $(2h+1)/(2f+1)$ a prime. If $M\cap \Oc\neq\emptyset$, then it is a unique semisimple conjugacy class in $M$ by Remark \ref{obs:ss-one-class}. By Remark \ref{obs:order-tori}, this class is represented in a torus of order coprime with $2^{2f+1}-1$, i.e., a torus as in the hypotheses. If $2f+1=1$, the intersection is non-trivial only if $|g|=5$ and by Remark \ref{obs:suz2} it consists of commuting elements. Arguing by induction on the number of prime factors of $2h+1$, we conclude that for any $K\leq \Gb$, the intersection $K\cap\Oc$ is either empty, a conjugacy class in $K$, or consists of commuting elements. \epf


\begin{obs}\label{obs:x1x2} Let $u\in\Gb$ be a non-trivial involution. Then, there exists $v\in\Oc_{u}^{\Gb}$ such that $(uv)^2\neq(vu)^2$. Indeed, since there is a unique conjugacy class of elements of order $2$ in $\Gb$,  we may always assume that $u\in {}^2\!B_2(2)$ and a direct computation using the realisation in Remark \ref{obs:suz2} shows that any $v\in\Oc_u^{^2\!B_2(2)}\setminus\{u\}$ has the required property. 
\end{obs}

\subsection{Nichols algebras over Suzuki groups}

In this subsection we consider Nichols algebras attached to simple Yetter-Drinfeld modules $M(\Oc,\rho)$ for $\Oc$ a kthulhu class in $\Gb$
and we prove Theorem \ref{thm:main} for simple Suzuki groups.

\begin{prop}\label{prop:suz_ss}
Let $g\in \Gb$ with $|g|\neq1$ and odd. Then $\dim\toba(\Oc_g,\rho)=\infty$ for every irreducible representation $\rho$ of $C_{\Gb}(g)$.
\end{prop}
\pf Every such element is semisimple, hence real by Remark \ref{obs:ss-real}. 
The claim follows from \cite[Lemma 2.2]{AZ}.
\epf

\begin{lema}\label{lem:acca}Assume $h>0$. Let $H=\langle \Tb, U(0,b),\,b\in\F_q\rangle= \Tb Z(U^-)$ and let $\Oc=\Oc_{U(0,1)}^H$. Then $\dim\toba(M(\Oc_{U(0,1)}^H,\rho))=\infty$  for every $\rho\in\Irr{C_H(U(0,1))}$.
\end{lema}
\pf The proof is obtained mutatis mutandis from the proof of  \cite[Proposition 3.1, case $c_2$]{FGV1}. We sketch it here for completeness' sake and we use notation and strategy from Subsection \ref{sec:abelian}. For $k\in\F_q^\times$ we consider the elements $t_k\trid U(0,1)=U(0,\varphi(k)^{-1})$. By Remark \ref{obs:kappa} all non-trivial involutions in $H$ are conjugate to $U(0,1)$ by an element in $\Tb$  and 
$A:=C_{H}(U(0,1))=Z(U^-)\simeq (\F_q, +)$  is abelian. We parametrise the elements in $\Oc=\Oc\cap A$ by elements in $H/Z(U^-)\simeq \Tb\simeq \F_q^\times$. Let $\chi$ be an irreducible representation of $Z(U^-)$, i.e.,  a group morphism $\chi\colon (\F_q,+)\to{\mathbb C}^\times$. The image is in $\{\pm1\}$.  Assume that $\dim\toba(M(\Oc,\chi))<\infty$. The coefficients of the  braiding associated with $(\Oc,\chi)$  are given by $q_{kl}=\chi(t_l^{-1}t_k\trid U(0,1))=\chi(U(0,\varphi(lk^{-1}))$ so $q_{kl}q_{lk}=\chi(U(0, \varphi(kl^{-1})+\varphi(lk^{-1}))$ for every $k,l\in\F_q^\times$. By \cite[Remark 1.1]{AZ} we necessarily have $q_{kk}=\chi(U(0,1))=-1$ for every $k\in\F_q^\times$.  Let  $k\in\F_q-\F_2$ and let $d$ be its multiplicative order.  Then $d\geq 3$ and if we had $q_{1k}q_{k1}=-1$, then we would have a cycle in the generalised Dynkin diagram associated with the braiding: $q_{1k}q_{k1}=q_{k k^2}q_{k^2k}=\cdots=q_{k^{d-1},1}q_{1 k^{d-1}}=-1$. This is excluded by \cite[Lemma 2.3]{FGV1}. Thus, $q_{1k}q_{k1}=\chi(U(0,\varphi(k)+\varphi(k^{-1})))=1$. However, the additive subgroup of $\F_q$ generated by the elements  of the form $\varphi(k)+\varphi(k^{-1})$, $k\in\F_q-\F_2$ contains $1$, so this would imply $\chi(U(0,1))=1$, a contradiction. 
\epf

\noindent\emph{Proof of Theorem \ref{thm:main} for simple Suzuki groups.} 
Proposition \ref{prop:type_Nichols} covers the cases of simple Yetter-Drinfeld modules associated with classes of elements of odd order or of order $4$  by Lemma \ref{lem:suz_order4} and Proposition \ref{prop:suz_ss}. We consider the class $\Oc$ of non-trivial involutions, represented by $U(0,1)$.  Lemma \ref{lem:acca} ensures that $\dim\toba(M(\Oc_{U(0,1)}^H,\chi))=\infty$  for every $\chi\in\Irr{C_H(U(0,1))}$. By \cite[Lemma 3.2]{AFGV-ampa},  $\dim\toba(M(\Oc,\rho))=\infty$ for every $\rho\in\Irr{C_{\Gb}(U(0,1))}$. We conclude by Remark \ref{obs:reduction}.\hfill$\Box$

\section{The Ree groups $^{2}\!F_4(2^{2h+1})$}

In this section $p=2$, $q=2^{2h+1}$, $h\geq0$, $\Gb={}^2\!F_4(q)$ a Ree group of type $F_4$. 

\subsection{Collapsing racks}
In this Subsection we list the kthulhu and non-kthulhu conjugacy classes in $\Gb$, when it is simple.  
We consider unipotent, semisimple and mixed classes separately.

\subsubsection{Unipotent classes}
We make use of the list of representatives of each conjugacy class in \cite[Table II]{Shi} and notation from \cite{ree2,AFGV-sporadic,FV} and \cite{HH}, with the understanding that our $q$ is $q^2$ therein. Roots in $\Phi^+$ are indicated as follows:  an index $j$ stands for $\varepsilon_j$; the symbol $i\pm j$ stands for $\varepsilon_i\pm\varepsilon_j$  and $1\pm2\pm3\pm4$ stands for $\frac{1}{2}(\varepsilon_1\pm\varepsilon_2\pm\varepsilon_3\pm\varepsilon_4)$. For $1\leq i\leq 12$, and $t\in \F_q$ there are elements $\alpha_i(t)\in U$ such that every $u \in U$ can be written uniquely as an ordered product $u=\prod_{i=1}^{12} \alpha_i(d_i)$ with $d_i\in\F_q$ for each $i=1,\ldots,12$. Commutation rules between  $\alpha_i(t)$ and $\alpha_j(t')$ are given in \cite{Shinoda2}: the case $(i,j)=(2,3)$ contained a mistake pointed out in \cite{Malle} and  corrected in \cite[Table 1]{HH}. We will also make use of the subsets $U_i :=\{\alpha_i(t)\mid t\in\F_q\}$ for $1\leq i\leq 12$ and of the subgroups $U_{\geq i}:=\prod_{j=i}^{12}U_j\unlhd U$.

\begin{lema}\label{lem:F_non_2A}Any non-trivial unipotent conjugacy class in $\Gb$, different from the one represented by $u_{1}$ in 
\cite[Table II]{Shi} is of type D. 
\end{lema}
\pf Each representative of the 19 unipotent conjugacy classes of $\Gb$  in \cite[Table II]{Shi} is defined over $\F_2$, i.e. it lies in the subgroup $^{2}\!F_4(2)$. Conjugacy classes in the Tits' group $^{2}\!F_4(2)'$ and in  $^{2}\!F_4(2)=\Aut(^{2}\!F_4(2)')$ are studied in \cite{AFGV-sporadic} and \cite{FV} respectively. By \cite[Table 2]{AFGV-sporadic} and \cite[Table 1]{FV} every conjugacy class of $^{2}\!F_4(2)$, apart from the one labeled by $2A$, is of type D. This class is represented by a non-trivial involution in $Z(U)$, hence it is the one represented by $u_{1}=\alpha_{12}(1)$. 
\epf

It has been shown in \cite[Proposition 4.1]{FV} that the class of $u_{1}$ is not of type D. Next Lemma deals with this class provided $h>0$.
\begin{lema}\label{lem:F_2A} Let $h>0$. The class $\Oc$ of $u_1$  is of type F.
\end{lema}
\pf Observe that $u_{1}=\alpha_{12}(1)=x_1(1)x_{1+2}(1)$. The  Weyl group element $s_{1-3}s_{2-4}$ lies in $C_W(\tau)$ so it has a representative $\dot{\sigma}$ in 
$^2\!F_4(2)\cap N_{\Gb}(\T)$. Hence $\dot{\sigma}\trid u_1=x_3(1)x_{3+4}(1)=\alpha_2(1)\in\Oc$. 
Let $\xi_j$ for $1\leq j\leq 4$ be distinct elements of $\F_q$. 
We consider the involutions $\alpha_3(\xi_j)=x_{1-2-3-4}(\xi_j^{2^h})x_{2-3}(\xi_j)\in U$ and we set
\begin{align*}
r_j:=\alpha_3(\xi_j)\trid u_1=\alpha_3(\xi_j)\trid\alpha_{2}(1)\in\Oc.
\end{align*}
Thus,
\begin{align*}
r_ir_j&=\alpha_3(\xi_i)\alpha_{2}(1)\alpha_3(\xi_i)\alpha_3(\xi_j)\alpha_{2}(1)\alpha_3(\xi_j)=\alpha_3(\xi_i)\alpha_{2}(1)\alpha_3(\xi_i+\xi_j)\alpha_{2}(1)\alpha_3(\xi_j)
\end{align*}
where we have used  Chevalley's commutator formula. 
Hence $r_ir_j=r_jr_i$ if and only if 
\begin{align*}
\alpha_3(\xi_i+\xi_j)\alpha_{2}(1)\alpha_3(\xi_i+\xi_j)\alpha_{2}(1)=\alpha_{2}(1)\alpha_3(\xi_i+\xi_j)\alpha_{2}(1)\alpha_3(\xi_j+\xi_i)
\end{align*}
and this happens if and only if the commutator of $\alpha_3(\xi_i+\xi_j)$ and $\alpha_2(1)$ is an involution. Making use of the commutation relations in \cite{Shinoda2,HH} we deduce
\begin{align*}
(\alpha_2(1)\alpha_3(\eta)\alpha_2(1)\alpha_3(\eta))^2=
(\alpha_5(\eta^{2^{h+1}})\alpha_6( \eta)\alpha_7(\eta)\alpha_8(\eta^{2^{h+1}+1})\alpha_9(\eta^{2^{h+1}+1}))^2\in\alpha_{10}(\eta)U_{\geq11}
\end{align*}
so $r_ir_j\neq r_jr_i$ whenever $i\neq j$. By direct computation:
\begin{align*}
r_i=\alpha_2(1)\alpha_5(\xi_i^{2^{h+1}})\alpha_6(\xi_i)\alpha_7(\xi_i)\alpha_8(\xi_i^{2^{h+1}+1})\alpha_9(u^{2^{h+1}+1})\in U_2U_5U_{\geq6}
\end{align*}
and $V=U_2U_5U_{\geq6}$ is a subgroup of $U$  with $V/U_{\geq 6}$ abelian.
Let $H=\langle r_1,\,r_2,\,r_3,\,r_4\rangle$. Then $H\leq V$ and so
$\Oc_{r_i}^H\subset \Oc_{r_i}^V\subset \alpha_2(1)\alpha_5(\xi_i^{2^{h+1}})U_{\geq6}$. Since $\xi_i^{2^{h+1}}=\xi_j^{2^{h+1}}$ only if $\xi_i=\xi_j$, the classes $\Oc_{r_i}^H$ for $i\neq j$ are disjoint and $\Oc$ is of type F.
\epf

\subsubsection{Mixed classes}

\begin{lema}\label{lem:mixed_no_inv}Let $\Oc$ be the class of an element $x\in \Gb$ with Jordan decomposition $x=x_sx_u$ with $x_s,x_u\neq1$. If $x_u^2\neq1$, then $\Oc$ is not kthulhu.
\end{lema}
\pf  Using strategy and notation from Remark \ref{obs:mixed} we see that $\Sigma$ can only be of type $A_1\times \tilde{A}_1$, $A_2\times \tilde{A}_2$, or $B_2$. By looking at the action of $F$ and the order of the centralisers given in \cite[Table IV]{Shi} we deduce that $ [C_{\G}(x_s),C_{\G}(x_s)]^F$ is either $\PSL_2(q)$, $\PSU_3(q)$, or $^2\!B_2(q)$ and $x_u$ lies in there. 
By  \cite[Proposition 5.1]{ACG-IV} and Lemma \ref{lem:suz_order4} we see that if $x_u^2\neq1$, then $\Oc_{x_u}^{[C_{\G}(x_s),C_{\G}(x_s)]^F}$ is not kthulhu. \epf

\begin{lema}\label{lem:ree_inv}Let $\Oc$ be the class of an element $x\in \Gb$ with Jordan decomposition $x=x_sx_u$ with $x_s,x_u\neq1$. Then $\Oc$ is not kthulhu.
\end{lema}
\pf By Lemma \ref{lem:mixed_no_inv} we may assume that $x_u^2=1$. We argue as in the proof of  \cite[Lemma 3.2]{ACG-V} to show that $\Oc$ is of type D. 
Since $w_0=-\id$, there is representative $\w_0\in N_{\Gb}(\tilde\T)$ for every $F$-stable maximal torus $\tilde\T$ containing $x_s$. Then,
$\w_0\trid x=x_s^{-1}x_u'$, where $x_s\neq x_s^{-1}$ because $x_s\neq1$ and $x_u'$ is a non-trivial involution in 
$K:=[C_{\G}(x_s^{-1}),C_{\G}(x_s^{-1})]^F=[C_{\G}(x_s),C_{\G}(x_s)]^F$. The latter is in turn isomorphic to
 $\PSL_2(q)$, $\PSU_3(q)$ or $^2\!B_2(q)$. All non-trivial involutions are conjugate in these groups, so $\Oc_{x_u'}^K=\Oc_{x_u}^K$. 
By  \cite[Lemma 3.6(a)]{ACG-IV}, \cite[Lemma 2.9]{ACG-V} and Remark \ref{obs:x1x2} there is $v\in \Oc_{x_u'}^K$ such that $(x_uv)^2\neq(vx_u)^2$.  
Thus, $s:=x_s^{-1}v\in\Oc_x^{\Gb}$ and we have: 
\begin{align*}(xs)^2\neq(sx)^2,&&\Oc_x^{\langle x,s\rangle}\subset x_s\Oc_{x_u}^K,&&\Oc_s^{\langle x,s\rangle}\subset x^{-1}_s\Oc_{x_u}^K, &&x_s\Oc_{x_u}^K\cap x^{-1}_s\Oc_{x_u}^K=\emptyset.
\end{align*}
The claim follows.\epf
 
\subsubsection{Semisimple classes}
 
The conjugacy classes of maximal tori in $\Gb$, the corresponding Weyl group elements, their orders and the order of their normalisers are listed in \cite[\S 3, Table III]{Shi}. They are represented by $\Tb_i$, for $1\leq i\leq 11$, with $|\Tb_i|=d_i$ as follows:
\begin{align*}
&d_1=(q-1)^2, &&d_2=(q-1)(q+1), \\
&d_{3,4}=(q-1)(q\pm\sqrt{2q}+1), &&d_5=(q^2+1),\\
& d_{6,7}=(q\pm \sqrt{2q}+1)^2, &&d_8=(q+1)^2,\\
 & d_9=(q^2-q+1),&&d_{10}=(q^2-\sqrt{2q^3}+q-\sqrt{2q}+1), \\
&  d_{11}=(q^2+\sqrt{2q^3}+q+\sqrt{2q}+1).
  \end{align*} 
We denote by $\T_i$, for $i\leq 11$ the corresponding $F$-stable maximal tori in $\G$.
Observe that
\begin{align*}
&(q-\sqrt{2q}+1)(q+\sqrt{2q}+1)=q^2+1,\\
& (q^3+1)=(q^2-q+1)(q+1),\\
&(q^2-\sqrt{2q^3}+q-\sqrt{2q}+1)(q^2+\sqrt{2q^3}+q+\sqrt{2q}+1)=q^4-q^2+1,\\
&(q^4-q^2+1)(q^2+1)=(q^6+1).
\end{align*}
According to \cite[\S 2.2]{Malle}, $\Gb$ contains  a unique conjugacy class of subgroups isomorphic to $^2\!B_2(q)\times {}^2\!B_2(q)$ and a unique conjugacy class of subgroups isomorphic to  $\SL_2(q)\times \SL_2(q)$. By looking at the maximal tori in $^2\!B_2(q)$ and $\SL_2(q)$ we see that
every torus $\Tb_i$ for $i\leq8$ is contained in a subgroup $M=M_1\times M_2$ with either $M_1\simeq M_2\simeq{} ^2\!B_2(q)$ or $M_1\simeq M_2\simeq \SL_2(q)$. This inclusion induces a decomposition of $\Tb_i$ into a  product of $2$ subtori $\Tb_{i,M_j}:=\Tb_i\cap M_j$ for $j=1,2$ whose orders follows from the decomposition of $d_i$ given above. If we write $x_s=(x_{1},x_{2})$ for an element in $\Tb_i$, we are referring to this decomposition. Also, we shall write $\T_{i,M_j}$ for the corresponding tori in $\Sp_4(\overline{\F_q})$ or $\SL_2(\overline{\F_q})$.

\begin{lema}\label{lem:order3}Let $\Oc$ be the class of an element of order $3$ in $\Gb$. Then, $\Oc$ is of type D.
\end{lema}
\pf According to \cite{Shi} there is a unique conjugacy class of elements of order $3$ in $\Gb$.  Recall that $|^2F'_4(2)|=2^{11}\cdot 3^3\cdot 5^2\cdot 13$, so the class $\Oc$ meets
 $^2F'_4(2)\leq \,^2\!F_4(2)\leq \,^2\!F_4(q)$. Since the only non-trivial class of $^2F'_4(2)$ which is not of type D consists of involutions \cite{AFGV-sporadic}, we have the statement.
\epf

\begin{lema}\label{lem:<9uno}Assume $h>0$. Let $\Oc=\Oc_{x_s}^{\Gb}$ for $x_s=(x_1,x_2)\in\Tb_i$ for $i\leq8$. If $x_1\neq1,\,x_2\neq1$, then, $\Oc$ is of type C.
\end{lema}
\pf   We consider the inclusion of $\Tb_i\leq M_1\times M_2$ with $M_j\simeq \SL_2(q)$ or $M_j\simeq \,^2\!B_2(q)$ for $j=1,2$.  The statement follows from Lemma \ref{lem:product} and Remark \ref{obs:ss-real}. 
\epf

\begin{lema}\label{lem:<9due}Assume $h>0$. Let $\Oc=\Oc_{x_s}^{\Gb}$ for $x_s=(x_1,x_2)\in\Tb_i\setminus 1$ for $i\leq8$. If $x_1=1$ or $x_2=1$, and $|x_s|\neq 3$, then, $\Oc$ is of type  C.
\end{lema}
\pf We assume that $x_1\neq1$, $x_2=1$, the other case is treated the same way.  If $\Oc_{x_s}^{\Gb}\cap \Tb_i$ contains $x_s'=(x_1',x_2')$ with $x_1'\neq1$, $x_2'\neq1$ then we invoke Lemma \ref{lem:<9uno}. 
Hence we assume from now on that $\Oc_{x_s}^{\Gb}\cap \Tb_i=(\Oc_{x_s}^{\Gb}\cap\Tb_{i,M_1})\cup(\Oc_{x_s}^{\Gb}\cap \Tb_{i,M_2})$.  The inclusion $x_s\in M_1\times M_2$  implies that $x_s\in \Tb_i$ for $i\in\{1,6,7,8\}$. Also, if 
$x_s=(x_1,1)\in \Tb_1$, then it lies in a split torus 
 in $\SL_2(q)=\PSL_2(q)$ or $^2\!B_2(q)$ and the claim follows either from \cite[Lemma 3.9]{ACG-III} or  Lemma \ref{lem:toro-split}.

 Thus, for the rest of the proof $x_s\in\Tb_i$ for  $i=6,7,8$, and $|x_s|\neq3$.  Observe that  $C_{\Gb}(x_s)\supset \Tb_{i,M_1}\times M_2$, so $C_{\Gb}(x_s)$ is not abelian. The structure of the centralisers of semisimple elements in $\Gb$ is described in  \cite[Theorem 3.2]{Shi} and most of them are tori. By inspection we see that  $x_s$ is necessarily conjugate to some $t_j$ from \cite[Table IV]{Shi}, with $j\in\{5,7,9\}$ and in these cases $C_{\Gb}(t_i)=\Tb_{i,M_1} M_2$  with $|\Tb_{i,M_1}|\in\{q+1,q\pm\sqrt{2q}+1\}$, so $x_1$ is regular in $M_1$. Observe that $M_1\simeq M_2\simeq \PSL_2(q)$ when $x_s\in\Tb_8$ and $x_s$ is conjugate to $t_5$ whereas  $M_1\simeq M_2\simeq\,^2\!B_2(q)$ when $x_s\in \Tb_6$ or $\Tb_7$ and in these cases $x_s$ is conjugate to $t_9$ or $t_7$. 

The inclusion $\Oc_{x_s}^{N_{\Gb}(\T_i)}\subseteq \Oc_{x_s}^{\Gb}\cap\Tb_i$
yields the inequality
\begin{align*}
|\Oc_{x_s}^{\Gb}\cap\Tb_i|&\geq |N_{\Gb}(\T_i)/(N_{\Gb}(\T_i)\cap C_{\Gb}(x_s))|=|N_{\Gb}(\T_i)/\Tb_{i,M_1}N_{M_2}(\T_{i,M_2})|\\
&=|N_{\Gb}(\T_i)/\Tb_i|/|N_{M_2}(\T_{i,M_2})/\Tb_{i,M_2}|.
\end{align*}
By \cite[\S 3, Table III]{Shi} the quotient $N_{\Gb}(\T_i)/\Tb_i$ has order $96$ (for $i=6,7$) or $48$ (for $i=8$). In addition, $|N_{M_2}(\T_{i,M_2})/\Tb_{i,M_2}|$ equals either $4$ (for $i=6,7$) or $2$ (for $i=8$). In all cases, $|\Oc_{x_s}^{\Gb}\cap\Tb_i|\geq 24$.

On the other hand, $|\Oc_{x_s}^{M_1\times M_2}\cap\Tb_i|=|\Oc_{x_s}^{M_1}\cap\Tb_{i,M_1}|$ cannot exceed the order of the Weyl group of $\Sp_{4}(\overline{\F_q})$ for $i=6,7$ and of $\SL_2(\overline{\F_q})$ for $i=8$, so $|\Oc_{x_s}^{M_1\times M_2}\cap\Tb_i|\leq 8$.
This shows that in our situation:
\begin{align*}\Oc_{x_s}^{M_1\times M_2}\cap \Tb_{i,M_1}=\Oc_{x_s}^{M_1\times M_2}\cap \Tb_i\subsetneq\Oc_{x_s}^{\Gb}\cap \Tb_i=(\Oc_{x_s}^{\Gb}\cap\Tb_{i,M_1})\cup(\Oc_{x_s}^{\Gb}\cap \Tb_{i,M_2}).
\end{align*}  The estimate $ |\Oc_{x_s}^{M_1\times M_2}\cap\Tb_{i}|\leq 8$ and its proof hold as well as if we replace $x_s$ by any $x_s'$ in $\Oc_{x_s}^{\Gb}\cap \Tb_{i,M_1}$ or in $\Oc_{x_s}^{\Gb}\cap \Tb_{i,M_2}$. Therefore the elements in $\Oc_{x_s}^{\Gb}\cap \Tb_i$ lie in at least $3$ distinct $(M_1\times M_2)$-orbits and each of these is contained either in $M_1$ or in $M_2$.  Without loss of generality we may assume that two of them are contained in $\Oc_{x_s}^{\Gb}\cap M_1$, say $\Oc_{s}^{M_1}$ and $\Oc_{s'}^{M_1}$, with $s,\,s'\in\Tb_{i,M_1}$ and by construction both regular in $M_1$. By Remark  \ref{obs:conditions} \eqref{item:andr} there is $g\in M_1$ such that $r:=g\trid s'\not\in C_{M_1}(s')=\Tb_{i,M_1}=C_{M_1}(s)$ so $[r,\,s]\neq1$. By Remark  \ref{obs:conditions} \eqref{item:bigger} the class is of type C with $H=\langle r,\,s\rangle$.
\epf

\begin{lema}\label{lem:9}Assume $h>0$. Let $\Oc=\Oc_{x_s}^{\Gb}$ for $x_s\in\Tb_9\setminus 1$ and $|x_s|\neq 3$.
Then  $\Oc$ is of type C.
\end{lema}
\pf According to  \cite{Malle}, $\Gb$ contains a subgroup isomorphic to $\SU_3(q)$, which contains a maximal torus of order $d_9$, so we may assume $x_s\in \Tb_9\leq \SU_3(q)$.
This torus consists of elements in $\SU_3(q)$ which are conjugate to diagonal matrices in  $\SL_3(\overline{\F_q})$  of the form $\textrm{diag}(x,x^{q^2},x^{-q})$, for $x^{q^2-q+1}=1$.
An element of this form could be real in $\SU_3(q)$ only if the set of its eigenvalues would coincide with its inverse set, which is impossible in our case,
so $\Oc_{x_s}^{\SU_3(q)}$ is not real. However, $\Oc$ is real in  $\Gb$ by Remark \ref{obs:ss-real}. Since $|x_s|\neq3$, it is not central in $\SU_3(q)$ so by Remark \ref{obs:conditions} \eqref{item:andr}, 
there is $g\in \SU_3(q)$ such that $[g\trid x_s^{-1},x_s]\neq 1$.
We take $r=g\trid x_s^{-1}$, $s=x_s$ and $H=\langle r,s\rangle\leq \SU_3(q)$ so $\Oc_s^H\neq\Oc_r^H$ and since $|x_s|$ is odd, $\Oc$ is of type C by Remark \ref{obs:conditions} \eqref{item:bigger}.    
\epf

\begin{lema}\label{lem:10-11}Let $\Oc=\Oc_{x_s}^{\Gb}$ for $x_s\in \Tb_i\setminus 1$ and $i=10,11$. If $|x_s|\neq13$, then  $\Oc$ is kthulhu.
\end{lema}
\pf A list of conjugacy classes of maximal subgroups for in $\Gb$ was computed in \cite{Malle} (note that $q^2$ there is our $q$). The list is exhaustive for $q\neq8$, whilst for $q=8$, one should add $3$ conjugacy classes of maximal subgroups isomorphic to ${\rm PGL}(13)$,  see  \cite[Remark 4.11]{Craven}. We use Remark \ref{obs:inductive_kthulhu} and consider the intersection of $\Oc$ with all maximal subgroups of $\Gb$. If $h=0$, then $ |\Tb_{10} ||\Tb_{11}|=13$, so there is nothing to prove.  First of all, we observe that  if $h=1$, then $ |\Tb_{10} |=37$ and $|\Tb_{11}|=109$, which are coprime with  $|{\rm PGL}_2(13)|=2^3\cdot 3\cdot 7\cdot13$. We now consider  the maximal groups listed in \cite{Malle} for $h>0$. 
Using coprimality of $q^4-q^2+1$ with $q$, $(q^2\pm1)$ and  $(q^2-q+1)$, we verify that $\Oc$ can have non-empty intersection only with $N_{\Gb}(\Tb_i)=\Tb_i\rtimes C_{12}$, for $i=10,11$ and $^2\!F_4(q_0)$ where $q_0=2^{2f+1}$ and $(2h+1)/(2f+1)$ is prime. Since $q\equiv2(3)$ and even, $(q^4-q^2+1,12)=1$, whence $\Oc\cap N_{\Gb}(\Tb_i)\subset \Tb_i$ is a commuting set. Assume $\Oc\cap\, ^2\!F_4(q_0)\neq\emptyset$. Then, it is a unique semisimple class in $^2\!F_4(q_0)$ by Remark \ref{obs:ss-one-class} and
by Remark \ref{obs:order-tori}, it has empty intersection with the maximal tori $\Tb_i'$ of $^2\!F_4(q_0)$ for $i\neq 10,11$. 
Hence, we are in the same hypotheses as above, with $q_0<q$. 
We proceed inductively on the number of prime factors of $2h+1$. When $2f+1=1$, our assumptions imply that $\Oc\cap \, ^2\!F_2(2)=\emptyset$.
\epf

\begin{lema}\label{lem:order13}Let $\Oc=\Oc_{x_s}^{\Gb}$ for $|x_s|=13$. Then,  $\Oc$ is not kthulhu.
\end{lema}
\pf If $x_s$ lies in  $\Tb_i$ for $i\leq 9$, the result follows from Lemmata \ref{lem:<9uno},   \ref{lem:<9due} and \ref{lem:9}. Assume $x_s\in\Tb_{j}$ for $j=10,11$. The torus $\Tb_j$ is cyclic and has empty intersection with all maximal tori of different order, so $\langle x_s\rangle$ is the only subgroup of order $13$ in $\Tb_j$ and all subgroups of this order are conjugate to it. Therefore $\Oc\cap \,^2\!F_4(2)'\neq\emptyset$, as $13$ divides the order of $^2\!F_4(2)'$.
 By \cite[Theorem II]{AFGV-sporadic} 
the classes contained in  $\Oc\cap \,^2\!F_4(2)'$ are of type D, whence so is $\Oc$.
\epf
\subsection{Nichols algebras over the Ree groups of type $F_4$}We are now in a position to prove Theorem  \ref{thm:main} for $\Gb$.

\noindent\emph{Proof of Theorem \ref{thm:main} for simple Ree groups of type $F_4$.} Proposition \ref{prop:type_Nichols} covers the cases of simple Yetter-Drinfeld modules associated with  unipotent classes by Lemmata \ref{lem:F_non_2A} and \ref{lem:F_2A};  mixed classes by Lemmata \ref{lem:mixed_no_inv} and \ref{lem:ree_inv} and semisimple classes represented in $\Tb_i$, for $i\leq 9$ by Lemmata \ref{lem:order3},  \ref{lem:<9uno},  \ref{lem:<9due}, \ref{lem:9}. The simple Yetter-Drinfeld modules associated with   classes represented in $\Tb_{10}$ or $\Tb_{11}$ are covered by the combination of Remark \ref{obs:ss-real} and \cite[Lemma 2.2]{AZ}. We conclude by Remark \ref{obs:reduction}. \hfill$\Box$
\section{The Ree groups $^{2}\!G_2(3^{2h+1})$}
In this section $p=3$, $q=3^{2h+1}$, $h\geq0$, $\Gb={}^2\!G_2(q)$ a Ree group of type $G_2$.  We fix a basis $\{\alpha,\beta\}$ of $\Phi$ with $\alpha$ short. We recall the list of maximal subgroups of $\Gb$ up to conjugation from \cite[Theorem C]{kle}:
\begin{enumerate}
 \item $\B^F$;
 \item the centraliser of a non-trivial involution $\sigma$, isomorphic to $\langle\sigma\rangle\times\PSL_2(q)$ (for $h>0$);
 \item the normaliser of a subgroup isomorphic to $C_2\times C_2$ (for $h>0$);
 \item $N_{\Gb}(\Tb_{s_\alpha s_\beta s_\alpha})\simeq \Tb_{s_\alpha s_\beta s_\alpha}\rtimes C_6$, 
 of order $6(q-\sqrt{3q}+1)$ (for $h>0$);
 \item $N_{\Gb}(\Tb_{s_\alpha s_\beta s_\alpha s_\beta s_\alpha})\simeq \Tb_{s_\alpha s_\beta s_\alpha s_\beta s_\alpha}\rtimes C_6$, of order $6(q+\sqrt{3q}+1);$
 \item $^2\!G_2(3^{2f+1})$ for $(2h+1)/(2f+1)$ a prime. 
\end{enumerate}
 If $h=0$, then $^2\!G_2(3)\simeq\PSL_2(8)\rtimes\langle \varphi\rangle$ where $\varphi$ acts as $\Fr_2$ on $\PSL_2(8)$ and, up to conjugation,  we have the additional maximal subgroups $\PSL_2(8)$ and the normaliser of the Sylow $2$-subgroup of upper triangular matrices in $\PSL_2(8)$, whose order is $2^3\cdot 3\cdot 7$.
 \begin{obs}\label{obs:maximalG2}We collect some properties of maximal subgroups of $\Gb$ and fix some notation: 
 \begin{enumerate}
 \item $\sigma$ will denote the involution $\sigma=\alpha^\vee(-1)\beta^\vee(-1)$, whose centraliser is the maximal subgroup $C_{\Gb}(\sigma)\simeq\langle\sigma\rangle\times\PSL_2(q)$. There is only one class of non-trivial involutions in $\G$  so by Remark \ref{obs:ss-one-class} there is only one class of non-trivial involutions in $\Gb$.
\item\label{item:semidir}  Assume $h>0$.  We provide a description of the normaliser of a subgroup isomorphic to $C_2\times C_2$ alternative to the one in \cite{kle}.
 There is a unique conjugacy class of such subgroups, \cite[p. 69]{W}. A representative is given by $\langle\sigma, \sigma'\rangle$ where $\sigma'$ the unique non-trivial involution in the (cyclic) maximal torus $\Tb'$ of $\PSL_2(q)\leq C_{\Gb}(\sigma)$ of order $\frac{q+1}{2}$.  The subgroup $\langle \sigma \rangle\times \Tb'$ is a maximal torus of $\Gb$ of order $q+1$ and there is only one class of tori of this order in $\Gb$, namely the one represented by $\Tb_{s_\alpha}$. We  set $\Tb_{s_\alpha}= \langle \sigma \rangle\times \Tb'$. 
By construction
$\Tb_{s_\alpha}\leq C_{\Gb}(\langle \sigma,\sigma'\rangle)\leq N_{\Gb}(\langle\sigma, \sigma'\rangle)$.
By \cite[List C]{kle} we have $N_{\Gb}(\langle\sigma, \sigma'\rangle)\simeq (C_2\times D_{(q+1)/2})\rtimes C_3$, so $\Tb_{s_\alpha}$ is normal of index $6$ in $N_{\Gb}(\langle\sigma, \sigma'\rangle)$. Simplicity of $\Gb$ and maximality of  $N_{\Gb}(\langle\sigma, \sigma'\rangle)$ give $N_{\Gb}(\langle\sigma, \sigma'\rangle)=N_{\Gb}(\Tb_{s_\alpha})$. Also $N_{\Gb}(\T_{s_\alpha})\leq N_{\Gb}(\Tb_{s_\alpha})$ and since $s_\alpha\tau$ is a rotation on $E$ we have $|C_W(s_\alpha\tau)|=6$. Hence, by order reason, 
$N_{\Gb}(\T_{s_\alpha})= N_{\Gb}(\Tb_{s_\alpha})$ and there exists an element $\varrho$ of order $6$ in $N_{\Gb}(\langle\sigma, \sigma'\rangle)$  such that  $\varrho^3\trid t=t^{-1}$ for every $t\in\Tb_{s_\alpha}$. 
Therefore $\langle \varrho\rangle\cap\Tb_{s_\alpha}=1$ and comparing orders we have $N_{\Gb}(\langle\sigma, \sigma'\rangle)=\Tb_{s_\alpha}\rtimes\langle\varrho\rangle$.
 \item\label{item:non-centralizza} Let $\xi\in\F_8^\times$ with $\xi^3+\xi+1=0$.
  Then $H:=\langle\left(\begin{smallmatrix}1&\xi\\
 0&1\end{smallmatrix}\right),\left(\begin{smallmatrix}1&\xi^2\\
0&1\end{smallmatrix}\right)\rangle\leq\PSL_2(8)\leq \,^2\!G_2(3)$
 is a representative of the conjugacy class of subgroups isomorphic to $C_2\times C_2$. Its normaliser in $^2\!G_2(3)$
  is generated by $\varphi$ and the subgroup of upper triangular unipotent matrices in $\PSL_2(8)$. 
  Observe that $\varphi$ permutes cyclically the three  non-trivial elements in $H$. Since all subgroups isomorphic to $H$ are conjugate in $\Gb$, it follows that also $\varrho^2$ as in \eqref{item:semidir} permutes $\sigma,\sigma'$ and $\sigma\sigma'$ cyclically. 
   \end{enumerate}
\end{obs} 
 
\subsection{Collapsing racks}
In this Subsection we list the kthulhu and non-kthulhu conjugacy classes in $\Gb$, when it is simple.  
We consider unipotent, semisimple and mixed classes separately.

\subsubsection{Unipotent classes}

 We see from \cite[Table 22.1.5]{LS} that $\G$ has 5 non-trivial unipotent classes: the regular one $\Oc_{reg}$,  the subregular one $\Oc_{subreg}$, represented by $x_{\beta}(1)x_{3\alpha+\beta}(1)$,  the class labeled by $(\tilde{A_1})_3$, represented by $x_{2\alpha+\beta}(1)x_{3\alpha+2\beta}(1)$, and the classes $\tilde{A_1}$ and ${A_1}$, represented by $x_\alpha(1)$ and $x_\beta(1)$, \cite[Table 2, Example 4.3]{sim}. The last two classes are interchanged by $F$, hence they do not intersect $\Gb$. The dimensions of the remaining ones are all distinct, hence they are all $F$-stable.  Regular unipotent elements have order $9$, all others have order $3$. By \cite[Table 22.1.5]{LS} the component group of $C_{\G}(u)$ is cyclic of order $2$ if $u\in\Oc_{subreg}$  and trivial if $u\in \Oc_{(\tilde{A_1})_3}$. 
In the first case $C_{\G}(u)$ is parted into two $F$-conjugacy classes, so $\Oc_{subreg}\cap\Gb$ is the union of two $\Gb$-classes, whereas $\Oc_{(\tilde{A_1})_3}\cap\Gb$ is a single unipotent conjugacy class. In addition $\Oc_{(\tilde{A_1})_3}$ is the unique unipotent class of dimension $8$ in $\G$, hence it is real, and therefore $\Oc_{(\tilde{A_1})_3}\cap\Gb$ is again so. The classes of $\varphi^{\pm1}$ in $^2\!G_2(3)$ are not real, hence these elements lie in $\Oc_{subreg}$ and represent the two unipotent classes in $\Gb$ contained therein.

\begin{lema}\label{lem:uni-reg-g2}
Let $h>0$. If $\Oc\subset \Oc_{reg}$, then $\Oc$ is of type D.
\end{lema}
\pf This is \cite[Proposition 3.7]{ACG-II}.\epf

\begin{lema}\label{lem:uni-subreg-reeg2}If $\Oc\subset \Oc_{subreg}\cup  \Oc_{(\tilde{A_1})_3}$, then $\Oc$ is kthulhu. 
\end{lema}
\pf We apply the strategy described in Remark \ref{obs:inductive_kthulhu} and consider the intersection of $\Oc$ with a maximal subgroup $M$ of $\Gb$. Notice that if $\Oc_{subreg}\cap H\neq\emptyset$ for some $H\leq \Gb$, then
$|\Oc_{\varphi}^{\Gb}\cap H|=|\Oc_{\varphi^{-1}}^{\Gb}\cap H|\neq0$. Let $M=\B^F$. The inclusion \begin{equation*}x_\alpha(\overline{\F_q}^\times)x_\beta(\overline{\F_q}^\times)\U_{\alpha+\beta}\U_{2\alpha+\beta}\U_{3\alpha+\beta}\U_{3\alpha+2\beta}\subset \Oc_{reg}\end{equation*} and $F$-invariance imply that  $\Oc\cap M\subset\U_{\alpha+\beta}\U_{3\alpha+\beta}\U_{2\alpha+\beta}\U_{3\alpha+2\beta}$ and all elements therein commute, \cite[III.6]{eno}. 

Let $M=C_{\Gb}(\sigma)$. It follows from  \cite[Lemmata 3.2 and 3.3, Corollary 3.4(ii)]{PST}, that $\Oc\cap M\subset \PSL_2(q)$ is  empty if $\Oc\subset  \Oc_{(\tilde{A_1})_3}$ and a single conjugacy class in $\PSL_2(q)$ if $\Oc$ is one of the two $\Gb$-classes contained in $\Oc_{subreg}$. Furthermore, when it is non-empty, the intersection of  $\Oc$ with a subgroup of $\PSL_2(q)$ is either a single conjugacy class or consists of commuting elements \cite[Lemma 3.5]{ACG-I}.

Let $M=N_{\Gb}(\langle \sigma,\sigma'\rangle)=\Tb_{s_\alpha}\rtimes\langle \varrho\rangle$. Observe that $\varrho^2\in C_{\Gb}(\varrho^3)$, the centraliser of an involution, and all elements of order $3$ therein are not real.  Hence, $\varrho^2\in\Oc_{subreg}\cap M=(\Oc_{\varphi}^{\Gb}\cap M)\cup(\Oc_{\varphi^{-1}}^{\Gb}\cap M)$. All elements of order $3$ in $M$ lie either in $\Tb_{s_\alpha}\varrho^2$ or in $\Tb_{s_\alpha}\varrho^4$, so  $(\Oc_{subreg}\cup  \Oc_{(\tilde{A_1})_3})\cap M\subset\Tb_{s_\alpha}\varrho^2\cup \Tb_{s_\alpha}\varrho^4$. Setting $M_1:= \Tb_{s_\alpha}\rtimes\langle \varrho^2\rangle$ we have
\begin{align*}
\Oc_{subreg}\cap M=\Oc_{subreg}\cap M_1=(\Oc_{\varphi}^{\Gb}\cap M_1)\cup (\Oc_{\varphi^{-1}}^{\Gb}\cap M_1).
\end{align*}
%
Observe that  $(q+1)/4$ is odd so  $\Tb_{s_\alpha}\simeq C_2\times C_2\times C_{(q+1)/4}$ and $\langle \sigma,\,\sigma'\rangle$ and $C_{(q+1)/4}$ are characteristic in $\Tb_{s_\alpha}$. We claim that $C_{\Tb_{s_\alpha}}(\varrho^2)=1$. Indeed, if $\varrho^2 t=t \varrho^2$ for some $t\in \Tb_{s_\alpha}$, then $\varrho^2$ would commute with the components of $t$ in $C_2\times C_2$ and $C_{(q+1)/4}$. The first component is trivial by Remark \ref{obs:maximalG2}\eqref{item:non-centralizza}, whereas the second one is trivial by \cite[p. 75]{W}. By Remark \ref{obs:split} \eqref{item:mauro2}, up to interchanging $\varrho$ and $\varrho^{-1}$,  we have 
\begin{align*}
\Oc_{subreg}\cap M=\Oc_{subreg}\cap M_1=\Tb_{s_\alpha}\varrho^2\cup\Tb_{s_\alpha}\varrho^4=\Oc_{\varrho^2}^{M_1}\cup\Oc_{\varrho^4}^{M_1}=\Oc_{\varrho^2}^{M}\cup\Oc_{\varrho^4}^{M}.
\end{align*}
Hence
$\Oc_{\varphi^{\pm1}}^{\Gb}\cap{M}=\Oc_{\varrho^{\pm2}}^M$ and $\Oc_{(\tilde{A_1})_3}\cap{M}=\emptyset$. Let now $H\leq M$ be  such that $\Oc_{\varphi^{\pm1}}^{\Gb}\cap H\neq\emptyset$. Replacing if needed $H$ by an $M$-conjugate of $H$ containing $\varrho^2$ we apply Remark \ref{obs:split} \eqref{item:mauro3} to deduce that $\Oc_{\varphi^{\pm1}}^{\Gb}\cap H=\Oc_{\varrho^{\pm2}}^H$.

Let $w={s_\alpha s_\beta s_\alpha}$ or $w={s_\alpha s_\beta s_\alpha s_\beta s_\alpha}$ and let $M$ be $N_{\Gb}(\Tb_w)=\Tb_w\rtimes \langle g\rangle$ for some 
$g\in\Gb$  with $|g|=6$. This case is similar to the case of $M=C_{\Gb}(\langle\sigma,\sigma'\rangle)$, but simpler. Here we use  \cite[Theorem, part (4)]{W}  to show that $C_{\Tb_w}(g^2)=1$ and proceed as before.

Let $M=\,^2\!G_2(q^{2f+1})$. There are three conjugacy classes of elements of order $3$ in $M$: the real one, which is $M\cap \Oc_{(\tilde{A_1})_3}$, and the two non-real ones, that are $M\cap \Oc_{\varphi}^{\Gb}$ and $M\cap \Oc_{\varphi^2}^{\Gb}$, so each intersection is a single conjugacy class in $M$ and we proceed inductively on the number of prime factors of $2h+1$. 

Finally, assume $q=3$. Let $M=\PSL_2(8)$. The only class of elements of order $3$ in $M$ is real, hence $M\cap \Oc_{(\tilde{A_1})_3}$ is this class in $M$ and $M\cap \Oc_{subreg}=\emptyset$. The intersection of $\Oc_{(\tilde{A_1})_3}$  with any subgroup of $\PSL_2(8)$ is either empty, a single conjugacy class,  or consists of commuting elements, \cite[Proposition 4.2, Case 2]{ACG-III}.

Let $M$ be the normaliser of a Sylow $2$-subgroup in $^2\!G_2(3)\simeq \PSL_2(8)\rtimes \langle \varphi\rangle$. Setting  $B_1:= \{
\left(\begin{smallmatrix}
\alpha&x\\
0&\alpha^{-1}\end{smallmatrix}\right)\mid x\in \F_8, \alpha\in \F_8^\times\}$ we  take $M=B_1\rtimes \langle\varphi\rangle$. Clearly, $\Oc_{subreg}\cap M\neq\emptyset$ and since $(|B_1|,3)=1$, we have the inclusion $(\Oc_{subreg}\cup \Oc_{(\tilde{A_1})_3})\cap M\subset B_1\varphi\cup B_1\varphi^{-1}$. Now, $C_{B_1}(\varphi)=\langle\left(\begin{smallmatrix}
1&1\\
0&1\end{smallmatrix}\right)\rangle$ and $\Oc_{\varphi^{\pm1}}^M\subset B_1\varphi^{\pm1}$, so $|\Oc_{\varphi^{\pm1}}^{M}|=|\Oc_{\varphi^{\pm1}}^{B_1}|=|B_1|/2$. 
The same argument shows that the orbits of the elements $\left(\begin{smallmatrix}
1&1\\
0&1\end{smallmatrix}\right)\varphi^{\pm1}$, whose order is $6$, have $|B_1|/2$ elements and lie in $B_1\varphi^{\pm1}$. Hence, 
$\Oc_{\varphi^{\pm1}}^{M}=\Oc_{\varphi^{\pm1}}^{\Gb}\cap M$ and $\Oc_{(\tilde{A_1})_3}\cap M=\emptyset$. Let now $H\leq M$ be  such that $H\cap\Oc_{subreg}\neq\emptyset$. Conjugating in $M$ we can always make sure that $\varphi\in H$, so $H=(B_1\cap H)\rtimes \langle \varphi\rangle$.  If  $\left(\begin{smallmatrix}
1&1\\
0&1\end{smallmatrix}\right)\not\in H$, then Remark \ref{obs:split} \eqref{item:mauro2} shows that $\Oc_{\varphi^{\pm1}}^H=(B_1\cap H)\varphi^{\pm1}=\Oc_{\varphi^{\pm1}}^{\Gb}\cap H$. If instead $\left(\begin{smallmatrix}
1&1\\
0&1\end{smallmatrix}\right)\in H$, then we use a counting argument as above to see that $\Oc_{\varphi^{\pm1}}^H=\Oc_{\varphi^{\pm1}}^{\Gb}\cap H$.
\epf

\subsubsection{Mixed classes}

\begin{lema}\label{lem:mixed_reeG}Let $h>0$ and let $\Oc=\Oc_x^{\Gb}$ where $x$ has Jordan decomposition $x=x_sx_u$ with $x_s,x_u\neq1$. Then $\Oc$ is of type D. 
\end{lema}
\pf Arguing as in Remark \ref{obs:mixed} we see that $[C_{\G}(x_s),C_{\G}(x_s)]$ is necessarily of type $A_1\times\tilde{A}_1$. In this case, $x_s$ is a non-trivial involution and we take $x_s=\sigma=\alpha^\vee(-1)\beta^\vee(-1)$. Thus $C_{\G}(\sigma)=\langle \T, \U_{\pm(\alpha+\beta)},\U_{\pm(3\alpha+\beta)}\rangle$, the roots $\alpha+\beta$ and $3\alpha+\beta$ are interchanged by $\theta$ and $C_{\Gb}(\sigma)\simeq\langle \sigma\rangle \times \PSL_2(q)$. Then $x_u$ can be chosen to be $x_{\alpha+\beta}(\epsilon)x_{3\alpha+\beta}(\epsilon)$ with $\epsilon=\pm1$, and the two choices represent two distinct conjugacy classes of elements of order $6$ in $\Gb$. By construction, there are no others.
Let $\F_q^\times=\langle\zeta\rangle$. We consider the elements $v,t,r,s$ and the subgroup $H$ as follows:
\begin{align*}
v&:=x_{2\alpha+\beta}(1)x_{3\alpha+2\beta}(1)\in \Gb\cap C_{\G}(\U_{\alpha+\beta}\U_{3\alpha+\beta}),\\
t&:=\alpha^\vee(\zeta^{3^h})\beta^\vee(\zeta)\in\Gb,\\
r&:=v\trid x_sx_u=\alpha^\vee(-1)\beta^\vee(-1)x_{\alpha+\beta}(\epsilon)x_{3\alpha+\beta}(\epsilon)x_{2\alpha+\beta}(1)x_{3\alpha+2\beta}(1)\in\Oc,\\
s&:=t\trid x_sx_u=\alpha^\vee(-1)\beta^\vee(-1)x_{\alpha+\beta}(\epsilon \zeta^{1-3^h})x_{3\alpha+\beta}(\epsilon\zeta^{3^{h+1}-1})\in\Oc,\\
H&:=\langle r,s\rangle\subset \langle \U_{\alpha+\beta}\U_{3\alpha+\beta}\U_{2\alpha+\beta}\U_{3\alpha+2\beta},\,\alpha^\vee(-1)\beta^\vee(-1)\rangle.
\end{align*}
where we have used that in characteristic $3$ there hold: $[x_s,s]=1$, $[\U_{2\alpha+\beta}\U_{3\alpha+2\beta},x_u]=1$ and that $x_s\trid x_{2\alpha+\beta}(\xi)x_{3\alpha+2\beta}(\eta)=x_{2\alpha+\beta}(-\xi)x_{3\alpha+2\beta}(-\eta)$ for every $\xi,\eta\in\overline{\F_3}$. Hence,
\begin{align*}
\Oc_{s}^H&\subset \alpha^\vee(-1)\beta^\vee(-1)x_{\alpha+\beta}(\epsilon \zeta^{1-3^h})x_{3\alpha+\beta}(\epsilon\zeta^{3^{h+1}-1})\U_{2\alpha+\beta}\U_{3\alpha+2\beta},\\
\Oc_{r}^H&\subset \alpha^\vee(-1)\beta^\vee(-1)x_{\alpha+\beta}(\epsilon)x_{3\alpha+\beta}(\epsilon)\U_{2\alpha+\beta}\U_{3\alpha+2\beta}
 \end{align*}
so $\Oc_s^H\cap\Oc_r^H=\emptyset$. In addition
\begin{align*}
(rs)^2&=\left(x_{\alpha+\beta}(\epsilon(1+\zeta^{1-3^h}))x_{3\alpha+\beta}(\epsilon(1+\zeta^{3^{h+1}-1}))x_{2\alpha+\beta}(-1)x_{3\alpha+2\beta}(-1)\right)^2\\
&=x_{\alpha+\beta}(2\epsilon(1+\zeta^{1-3^h}))x_{3\alpha+\beta}(2\epsilon(1+\zeta^{3^{h+1}-1}))x_{2\alpha+\beta}(1)x_{3\alpha+2\beta}(1),\\
(sr)^2&=\left(x_{\alpha+\beta}(\epsilon(1+\zeta^{1-3^h}))x_{3\alpha+\beta}(\epsilon(1+\zeta^{3^{h+1}-1}))x_{2\alpha+\beta}(1)x_{3\alpha+2\beta}(1)\right)^2\\
&=x_{\alpha+\beta}(2\epsilon(1+\zeta^{1-3^h}))x_{3\alpha+\beta}(2\epsilon(1+\zeta^{3^{h+1}-1}))x_{2\alpha+\beta}(2)x_{3\alpha+2\beta}(2),
\end{align*}
so $(rs)^2\neq(sr)^2$  and $\Oc$ is of type D.
\epf

\subsubsection{Semisimple classes}

There are four $\Gb$-conjugacy classes of maximal tori represented by $\Tb$, $\Tb_{s_\alpha}$, $\Tb_{s_\alpha s_\beta s_\alpha}$ and $\Tb_{s_\alpha s_\beta s_\alpha s_\beta s_\alpha}$ of order $q\mp1$, $q\mp\sqrt{3q}+1$, respectively. Their orders are mutually coprime in all cases except from $(|\Tb|,|\Tb_{s_\alpha}|)=2$.
We realise $\Tb$ and $\Tb_{s_\alpha}$ in $C_{\Gb}(\sigma)$ as the direct product of $\langle \sigma\rangle$ and a maximal torus in $\PSL_2(q)$, so $\Tb\cap\Tb_{s_\alpha}$ is a cyclic group of order $2$. The tori $\Tb$, $\Tb_{s_\alpha s_\beta s_\alpha}$ and $\Tb_{s_\alpha s_\beta s_\alpha s_\beta s_\alpha}$  are  cyclic. 
%

\begin{obs}\label{obs:ss_conjugacy}Let $t\in \Gb$ be semisimple, and such that $t^2\neq1$. Then $[C_{\G}(t),C_{\G}(t)]$ is not of type $A_1\times\tilde{A}_1$, so it is trivial. In other words,  $t$ is regular in $\G$ and it lies in a unique maximal torus of $\G$. Hence, if $\Tb_w=\T_w^F$ is a maximal torus in $\Gb$ and $t_1,t_2\in\Tb_w$ satisfy $g\trid t_1=t_2$ for some $g\in \Gb$ and $t_1^2\neq1$, then $g\trid {\T}_w=g\trid C_{\G}(t_1)=C_{\G}(t_2)={\T}_w$. Hence, $g\in N_{\Gb}(\T_w)$ and 
$\left|\Oc_{t_1}^{\Gb}\cap\Tb_w\right|=\left|\Oc_{t_1}^{N_{\Gb}(\T_w)}\right|=\left|N_{\Gb}(\T_w)/\Tb_w\right|=|C_W(w\tau)|$.
\end{obs}
\begin{lema}\label{lem:ree-semis-T1}
Let $h>0$. If $x_s\in\Tb\setminus\{1\}$, then $\Oc_{x_s}^{\Gb}$ is of type D. 
\end{lema}
\pf There is only one class of non-trivial involutions in $\Gb$ so if $|x_s|=2$ its class is represented by $\sigma'\in\PSL_2(q)$. If, instead, $|x_s|>2$, we may assume that $\sigma\in \Tb$ so $\Tb\leq C_{\Gb}(\sigma)\simeq \langle \sigma\rangle\times\PSL_2(q)$ and $x_s$ is conjugate either to $y$ or to $\sigma y$ for some $y\neq 1$ in a maximal torus of order $(q-1)/2$ for $\PSL_2(q)$. By   \cite[Corollary 3.5, Lemma 3.9]{ACG-III}  the racks $\Oc_{\sigma'}^{\PSL_2(q)}$ and  $\Oc_{\sigma y}^{\langle\sigma\rangle\times \PSL_2(q)}\simeq \Oc_y^{\PSL_2(q)}$ are of type D. We conclude by using Proposition \ref{prop:type_Nichols} \eqref{item:inj-proj}.\epf

\begin{lema}\label{lem:order2}Let $h>0$. If $|x_s|=7$, then $\Oc_{x_s}^{\Gb}$ is of type D. 
\end{lema}
\pf 
There is precisely one maximal torus up to conjugacy containing elements of order $7$ and by the structure of the tori, it contains exactly one subgroup of order $7$. Thus, all subgroups of order $7$ are conjugate in $\Gb$. We recall that there is a subgroup isomorphic to $\PSL_2(8)$ in $^2\!G_2(3)\leq \Gb$. It contains a subgroup of order $7$, namely its split torus, which intersects $\Oc_{x_s}^{\Gb}$. We conclude by invoking \cite[Lemma 3.9]{ACG-III} and Proposition \ref{prop:type_Nichols} \eqref{item:inj-proj}.
\epf
\begin{lema}\label{lem:semis-T34-kthulhu}
Let $w=s_\alpha s_\beta s_\alpha$ or $s_\alpha s_\beta s_\alpha s_\beta s_\alpha$ and let $x_s\in\Tb_w\setminus\{1\}$.  If $|x_s|\neq 7$, then $\Oc_{x_s}^{\Gb}$ is kthulhu.  
\end{lema}
\pf We use the strategy from Remark \ref{obs:inductive_kthulhu}.  If $h=0$, then $ |\Tb_w|=1$ or $7$ , so there is nothing to prove.  We assume $h > 0$. The order of $x_s$ divides $q^2-q+1$ so it is odd and coprime with $q\pm1$ and $q$. Hence, the only maximal subgroups intersecting $\Oc_{x_s}^{\Gb}$ are $N_{\Gb}(\Tb_w)$ and $^2\!G_2(3^{2f+1})$ with $(2h+1)/(2f+1)$ a prime number. In the first case $\Oc_{x_s}^{\Gb}\cap M$ consists of commuting elements; in the second case it is 
a single conjugacy class by Remark \ref{obs:ss-one-class}. In addition,  Remark \ref{obs:order-tori} shows that $x_s$ cannot lie in a torus of order $3^{2f+1}\pm1$. 
We proceed inductively on the number of prime factors of $2h+1$. When $2f+1=1$, our assumptions imply that $\Oc\cap \, ^2\!G_2(3)=\emptyset$.
\epf

\begin{lema}\label{lem:semis_T2}Let $h>0$. If $x_s\in\Tb_{s_\alpha}\setminus\{1\}$ and $|x_s|\neq 2$, then $\Oc_{x_s}^{\Gb}$ is of type C.  
\end{lema}
\pf We recall from Remarks \ref{obs:maximalG2} \eqref{item:semidir}
and \ref{obs:ss_conjugacy} that $\Tb_{s_\alpha}=\langle\sigma\rangle\times\Tb'\leq\langle\sigma\rangle\times\PSL_2(q)$ and that $x_s$ is regular. 
Observe that $s_\alpha\tau$ is a rotation on $E$ so $C_W(s_\alpha\tau)$ is cyclic of order $6$. By Remark \ref{obs:ss_conjugacy} we have 
\begin{align*}\left|\Oc_{x_s}^{\Gb}\cap\Tb_{s_\alpha}\right|=6>2=\left|\{x_s,x_s^{-1}\}\right|=\left|\Oc_{x_s}^{\langle\sigma\rangle\times\PSL_2(q)}\cap\Tb_{s_\alpha}\right|.\end{align*}
Hence, there is $s\in (\Oc_{x_s}^{\Gb}\cap\Tb_{s_\alpha})\setminus\Oc_{x_s}^{\langle\sigma\rangle\times\PSL_2(q)}$. 
By Remark \ref{obs:conditions} \eqref{item:andr} there is $g\in\langle\sigma\rangle\times\PSL_2(q)$ such that $r:=g\trid x_s\not\in C_{\Gb}(x_s)=\Tb_{s_\alpha}=C_{\Gb}(s)$. If we set $H:=\langle r,\,s\rangle\leq\langle\sigma\rangle\times\PSL_2(q)$, then $\Oc_s^H\neq\Oc_r^H$. If $|x_s|$ is odd the class is of type C by
Remark \ref{obs:conditions} \eqref{item:bigger}. If $|x_s|$ is even, then we decompose $s=s_es_o$ into its $2$-part and $2$-regular part, so $s_o\in H$, all prime factors of its order are $\geq5$ and $s_o$ is regular in $\G$ by Remark \ref{obs:ss_conjugacy}. We thus have $\left|\Oc_{r}^H\right|\geq \left|\Oc_{r}^{\langle s_o\rangle}\right|\geq5$ and therefore $\Oc_{x_s}^{\Gb}$ is of type C.
\epf

\subsection{Nichols algebras over the Ree groups of type $G_2$}

In this subsection we consider Nichols algebras attached to simple Yetter-Drinfeld modules $M(\Oc,\rho)$ for $\Oc$ a kthulhu class in $\Gb$
and we prove Theorem \ref{thm:main} for simple Ree groups of type $G_2$.

\begin{prop}\label{prop:ree-g2-uni}
Assume $h>0$. Let $g\in \Gb$ be a non-trivial unipotent element. Then $\dim\toba(\Oc^{\Gb}_g,\rho)=\infty$ for every irreducible representation $\rho$ of $C_{\Gb}(g)$.
\end{prop}
\pf If $g\in \Oc_{reg}$, this follows from Proposition \ref{prop:type_Nichols} (1), Lemma \ref{lem:uni-reg-g2} and \cite[Theorem 3.6]{AFGV-ampa}. If $g\in\Oc_{(\tilde{A_1})_3}$, then $\Oc_g^{\Gb}$ is real and of odd order and the claim follows from \cite[Lemma 2.2]{AZ}. Assume now $g=\varphi\in \Oc_{\varphi}^{\Gb}$, the case of $g=\varphi^{-1}$ is treated similarly. 
We will show that $\dim\toba(\Oc_{\varphi}^{^2\!G_2(3)},\rho)=\infty$ for every irreducible representation $\rho$ of $C_{^2\!G_2(3)}(\varphi)$ and apply \cite[Lemma 3.2]{AFGV-ampa}. Recall that $\Gbtre=\PSL_2(8)\rtimes\langle\varphi\rangle$, so setting $N:=\PSL_2(8)$ we see that $\Oc_{\varphi^{\pm1}}^{^2\!G_2(3)}\subset N\varphi^{\pm1}$. The elements of order $3$ in $N$ are real, so the real class of elements of order $3$ in $\Gbtre$ is all contained therein. Thus, an element of order $3$ lies in 
$\Oc_{\varphi^{\pm1}}^{^2\!G_2(3)}$ if and only if it lies in $N\varphi^{\pm1}$. 
We proceed as outlined in Subsection \ref{sec:abelian}, from which we adopt notation. We consider 
\begin{align*}C_{\Gbtre}(\varphi)\simeq\SL_2(2)\times\langle\varphi\rangle\simeq\s_3\times \langle \varphi\rangle,&&A:={\mathbb A}_3\times\langle\varphi\rangle\leq C_{\Gbtre}(\varphi).\end{align*}
The intersection $\Oc_{\varphi}^{\Gbtre}\cap C_{\Gbtre}(\varphi)=\{\varphi,(123)\varphi,(132)\varphi\}=\Oc_{\varphi}^{\Gbtre}\cap A$ is a commuting set and we put $x_0=\varphi$, $x_1=(123)\varphi$, $x_2=(132)\varphi$, so $x_0x_1x_2=1$. 
We claim that there is $g\in\Gbtre$ such that $g^j\trid x_i=x_{i+j\mod 3}$ for $j\in{\mathbb Z}$, $i=0,1,2$. 
Let $z\in\Gbtre$ be such that $z\trid x_0=x_1$, so $[z\trid x_0,x_0]=1$, whence $z^{-1}\trid x_0\in\Oc_{\varphi}^{\Gbtre}\cap C_{\Gbtre}(\varphi)$ and $z^{-1}\trid x_0\neq x_0$. If $z^{-1}\trid x_0=x_2$, then 
$z\trid x_2=x_0$ and so $z\trid x_1=z\trid (x_0x_2)^{-1}=(x_1x_0)^{-1}=x_2$ and we put $g=z$. If, instead, $z^{-1}\trid x_0=x_1$, then $z\trid x_2=x_2$  
and we consider $y\in\Gbtre$ such that $y\trid x_0=x_2$ and repeat the argument. Then either $g=y$ will do, or $g=yz$ will do. 
Let $\rho$ be an irreducible representation of $C_{\Gbtre}(x_0)$ and let  ${\mathbb C}v$ be any line stabilised by $A$. Let $\rho(x_i)v=\omega_iv$ for $i=0,1,2$. 
We have $\omega_i^3=1$ and $\omega_0\omega_1\omega_2=1$. By \eqref{eq:braiding-abelian} the braided vector subspace $M_A={\rm span}_{\mathbb C}\{g^i\otimes v,\,i=0,1,2\}$ of $M(\Oc_{\varphi}^{\Gbtre},\rho)$ has braiding
\begin{align*}
c((g^i\otimes v)\otimes (g^j\otimes v))=q_{ij}(g^j\otimes v)\otimes (g^i\otimes v),&&
q_{ij}v=\rho(g^{-j+i}\trid x_0)v=\omega_{i-j{\rm mod}3}v.
\end{align*}
This gives
\begin{align*}
q_{01}q_{10}=q_{12}q_{21}=q_{02}q_{20}=\omega_1\omega_2=\omega_0^2,&&\omega_{ii}=\omega_0, i=0,1,2.
\end{align*}
If $\omega_0=1$, then $\dim\toba(M(\Oc_{\varphi}^{\Gbtre},\rho))=\infty$ by \cite[Remark 1.1]{AZ}. If, instead, $\omega_0$ is a primitive third root of $1$, then the generalized Dynkin diagram of $M_A$ is connected and does not occur in \cite[Table 2]{H}. This implies $\dim\toba(M_A)=\infty$, whence again $\dim\toba(M(\Oc^{\Gbtre}_{\varphi},\rho))=\infty$.
\epf

\noindent\emph{Proof of Theorem \ref{thm:main} for simple Ree groups of type $G_2$.} Proposition \ref{prop:type_Nichols} covers the cases of simple Yetter-Drinfeld modules associated with  unipotent classes by Proposition \ref{prop:ree-g2-uni}; 
 mixed classes by Lemma \ref{lem:mixed_reeG} and semisimple classes represented in  $\Tb_w$ for $w=1$ or $s_\alpha$ by Lemmata \ref{lem:ree-semis-T1}, \ref{lem:order2} and \ref{lem:semis_T2}.
 The simple Yetter-Drinfeld modules associated with   classes represented in  $\Tb_w$ for $w=s_\alpha s_\beta s_\alpha$ or $s_\alpha s_\beta s_\alpha s_\beta s_\alpha$ are covered by the combination of Remark \ref{obs:ss-real} and \cite[Lemma 2.2]{AZ}. We conclude by Remark \ref{obs:reduction}.
\hfill$\Box$

\begin{obs}
The non-simple groups $^2\!F_4(2)$ and $^2\!G_2(3)$ have simple derived subgroup $\Gb'$. The group $^2\!F_4(2)'$ has been dealt with in \cite{AFGV-sporadic,FV}, whereas $^2\!G_2(3)'\simeq\PSL_2(8)$ has been treated in \cite{FGV1}. In both cases, $\dim\toba(V)=\infty$ for every $V\in\,^{\Gb'}_{\Gb'}\!\mathcal{YD}$.
The group $^2\!B_2(2)$ has order $20$ and all its classes are kthulhu, see Remark \ref{obs:suz2}. Its derived subgroup is cyclic of order $5$, so it is simple and abelian, and there do exist finite-dimensional pointed Hopf algebras over $C_5$, for example the Taft algebras in \cite{Taft}, see also \cite[Theorem 1.3]{AS}.
\end{obs}


\section{Acknowledgements}
The authors could benefit from several discussions with N. Andruskiewitsch and G. A. Garc\'ia, and are grateful to them. In particular, the content of Remark \ref{obs:conditions} (1) was pointed out to G.C. by N. Andruskiewitsch.


\begin{thebibliography}{29}


\bibitem{ACG-I}  N. Andruskiewitsch, G. Carnovale, G. A. Garc\'ia.
\emph{Finite-dimensional pointed Hopf algebras over finite simple groups of Lie type  I. Unipotent classes in  $\PSL_n(q)$},
 J. Algebra, 442, 36--65 (2015).

\bibitem{ACG-II}  \bysame,\emph{Finite-dimensional pointed Hopf algebras over finite simple groups of Lie type  II. 
Unipotent classes in  symplectic groups}, Commun. Contemp. Math. 18, No. 4, Article ID 1550053, 35 pp. (2016).

\bibitem{ACG-III}  \bysame, \emph{Finite-dimensional pointed Hopf algebras over finite simple groups of Lie type  III. 
Semisimple classes in $\PSL_n(q)$}, Rev. Mat. Iberoam.  \textbf{33(3)}, 995--1024, (2017).

\bibitem{ACG-IV} \bysame, \emph{Finite-dimensional pointed Hopf algebras over finite simple groups of Lie type  IV. 
Unipotent classes in Chevalley and Steinberg groups}, Algebras and Representation Theory, https://doi.org/10.1007/s10468-019-09868-6

\bibitem{ACG-V} \bysame, \emph{Finite-dimensional pointed Hopf algebras over finite simple groups of Lie type  V. Mixed classes in Chevalley and Steinberg groups}, to appear in Manuscripta Mathematica.


\bibitem{AFGV-ampa} N. Andruskiewitsch, F. Fantino, M.
Gra\~na,  L. Vendramin, \emph{Finite-di\-mensional pointed
Hopf algebras with alternating groups are trivial},
Ann. Mat. Pura Appl. (4),  \textbf{190}  (2011), 225--245.

\bibitem{AFGV-sporadic} \bysame, \emph{Pointed Hopf algebras over the sporadic simple groups}. J.
Algebra \textbf{325} (2011), 305--320.

\bibitem{AG}N. Andruskiewitsch, M. Gra\~na, {\emph From racks to pointed Hopf algebras}, Adv. Math. 178 (2003), 177--243. 

\bibitem{AS}N. Andruskiewitsch,  H. Schneider,
\emph{Finite quantum groups and Cartan matrices}, 
Adv. Math. 154 (2000), 1--45.


\bibitem{AZ}  N. Andruskiewitsch,  S. Zhang,
\emph{On pointed Hopf algebras associated to some conjugacy
	classes in $\mathbb S_n$}, Proc. Amer. Math. Soc. \textbf{135}  (2007),
2723--2731.

\bibitem{Carter1} R.W. Carter,
\emph{Simple Groups of Lie Type, Reprint of the 1972 original},
John Wiley (1989).



\bibitem{Craven}D.A. Craven, \emph{The maximal subgroups of the exceptional groups $F_4(q)$, $E_6(q)$ and ${}^2\!E_6(q)$ and related almost simple groups}, Invent. Math. \textbf{234}  (2023), 637--719.

\bibitem{eno}{H. Enomoto},
\emph{The conjugacy classes of Chevalley groups of type $G_2$ over finite fields of characteristic $2$ or $3$,}
J. Fac. Sci. Univ. Tokyo, Sect. I,  \textbf{16}  (1970), 497-512.

\bibitem{FV} F. Fantino, L.  Vendramin, On twisted conjugacy classes of type D in sporadic simple groups,
\emph{Contemp. Math.} \textbf{585} (2013), 247--259.


\bibitem{FGV1} S. Freyre, M. Gra\~na and L. Vendramin,
\emph{On Nichols algebras over $\SL(2,q)$ and $\GL(2,q)$}. J. Math. Phys. \textbf{48}, (2007) 123513.


\bibitem{H}I. Heckenberger,
\textit{Classification of arithmetic root systems of rank 3},
in: Actas del ``XVI Coloquio Latinoamericano
de \'Algebra'', Colonia, Uruguay, 2005, 227--252.

\bibitem{HH} F. Himstedt, S.-C. Huang,
\emph{Character table of a Borel subgroup of the Ree groups $^{2}\!F_4(q^2)$},
LMS J. Comput. Math.  \textbf{12} (2009) 1--53

\bibitem{kle}P. B. Kleidman,
\emph{The maximal subgroups of the Chevalley groups $G_2(q)$ with $q$ odd, the Ree groups $^2\!G_2(q)$, and their automorphism groups},
J. Algebra 117, (1988),  30--71.

\bibitem{LS} M. W. Liebeck, G. Seitz,
\emph{Unipotent and Nilpotent Classes in Simple Algebraic Groups and Lie Algebras}, SURV 180, AMS, 2012.

\bibitem{Malle}G. Malle,
\emph{The maximal subgroups of $^2\!F_4(q^2)$},
J. Algebra 139, 52--69, (1991).

\bibitem{MT} G. Malle and D. Testerman,
\emph{Linear Algebraic Groups and Finite Groups of Lie Type},
Cambridge Studies in Advanced Mathematics \textbf{133} (2011).

\bibitem{PST}R. Proud, J. Saxl, D. Testerman,
\emph{Subgroups of type $A_1$ containing a fixed unipotent element in an algebraic group},
J. Algebra 231, 53--66, (2000).

\bibitem{ree1}R. Ree,
\emph{A family of simple groups associated with the simple Lie algebra of type (G2)}, Amer. J. Math. 83, 432--462 (1961). 

\bibitem{ree2}\bysame,
\emph{A family of simple groups associated with the simple Lie algebra of type (F4)}, Amer. J. Math. 83,  401--420, (1961). 

\bibitem{Shi}K. Shinoda,
\emph{The conjugacy classes of the finite Ree groups of type $(F_4)$},
J. Fac. Sci. Univ. Tokyo 22, (1975),  1--15.

\bibitem{Shinoda2}K. Shinoda,
\emph{A characterization of odd order extensions of the {R}ee groups
              $^{2}F_{4}(q)$},
J. Fac. Sci. Univ. Tokyo Sect. I A Math. 22, (1975),  79--102.

\bibitem{sim}I.I. Simion,
\emph{Double centralizers of unipotent elements in simple algebraic groups of type $G_2$, $F_4$ and $E_6$},
J. Algebra 382 (2013), 335--367. 


\bibitem{sp-st} T. A. Springer, R. Steinberg,
\emph{Conjugacy Classes, Seminar on Algebraic Groups and Related Finite Groups},
167--266, LNM 131, Springer,  (1970).

\bibitem{Su} M. Suzuki,
\emph{On a class of doubly transitive groups}.
Ann. of Math. (2) 75 (1962), 105--145.

\bibitem{Taft}E. Taft,
\emph{The order of the antipode of finite-dimensional Hopf algebra}, Proc. Natl. Acad. Sci. U.S.A. 68 (1971) ,2631--2633.

\bibitem{W} H.N. Ward,
\emph{On Ree's series of simple groups}.
Trans. Amer. Math. Soc. 121 (1966), 62--89.



\end{thebibliography}
\end{document}